
\ifx\shlhetal\undefinedcontrolsequence\let\shlhetal\relax\fi
\def\fmtname{AmS-TeX}

\def\fmtversion{2.1}
\catcode`\@=11
\ifx\amstexloaded@\relax\catcode`\@=\active
  \endinput\else\let\amstexloaded@\relax\fi
\newlinechar=`\^^J
\def\W@{\immediate\write\sixt@@n}
\def\CR@{\W@{^^J\fmtname - Version \fmtversion^^J%
COPYRIGHT 1985, 1990, 1991 - AMERICAN MATHEMATICAL SOCIETY^^J%
Use of this macro package is not restricted provided^^J%
each use is acknowledged upon publication.^^J}}
\CR@ \everyjob{\CR@}
\message{Loading definitions for}
\message{misc utility macros,}
\toksdef\toks@@=2
\long\def\rightappend@#1\to#2{\toks@{\\{#1}}\toks@@
 =\expandafter{#2}\xdef#2{\the\toks@@\the\toks@}\toks@{}\toks@@{}}
\def\alloclist@{}
\newif\ifalloc@
\def\showallocations{{\def\\{\immediate\write\m@ne}\alloclist@}\alloc@true}
\def\alloc@#1#2#3#4#5{\global\advance\count1#1by\@ne
 \ch@ck#1#4#2\allocationnumber=\count1#1
 \global#3#5=\allocationnumber
 \edef\next@{\string#5=\string#2\the\allocationnumber}%
 \expandafter\rightappend@\next@\to\alloclist@}
\newcount\count@@
\newcount\count@@@
\def\FN@{\futurelet\next}
\def\DN@{\def\next@}
\def\DNii@{\def\nextii@}
\def\RIfM@{\relax\ifmmode}
\def\RIfMIfI@{\relax\ifmmode\ifinner}
\def\setboxz@h{\setbox\z@\hbox}
\def\wdz@{\wd\z@}
\def\boxz@{\box\z@}
\def\setbox@ne{\setbox\@ne}
\def\wd@ne{\wd\@ne}
\def\iterate{\body\expandafter\iterate\else\fi}
\def\err@#1{\errmessage{AmS-TeX error: #1}}
\newhelp\defaulthelp@{Sorry, I already gave what help I could...^^J
Maybe you should try asking a human?^^J
An error might have occurred before I noticed any problems.^^J
``If all else fails, read the instructions.''}
\def\Err@{\errhelp\defaulthelp@\err@}
\def\eat@#1{}
\def\in@#1#2{\def\in@@##1#1##2##3\in@@{\ifx\in@##2\in@false\else\in@true\fi}%
 \in@@#2#1\in@\in@@}
\newif\ifin@
\def\space@.{\futurelet\space@\relax}
\space@. %
\newhelp\athelp@
{Only certain combinations beginning with @ make sense to me.^^J
Perhaps you wanted \string\@\space for a printed @?^^J
I've ignored the character or group after @.}
{\catcode`\~=\active 
 \lccode`\~=`\@ \lowercase{\gdef~{\FN@\at@}}}
\def\at@{\let\next@\at@@
 \ifcat\noexpand\next a\else\ifcat\noexpand\next0\else
 \ifcat\noexpand\next\relax\else
   \let\next\at@@@\fi\fi\fi
 \next@}
\def\at@@#1{\expandafter
 \ifx\csname\space @\string#1\endcsname\relax
  \expandafter\at@@@ \else
  \csname\space @\string#1\expandafter\endcsname\fi}
\def\at@@@#1{\errhelp\athelp@ \err@{\Invalid@@ @}}
\def\atdef@#1{\expandafter\def\csname\space @\string#1\endcsname}
\newhelp\defahelp@{If you typed \string\define\space cs instead of
\string\define\string\cs\space^^J
I've substituted an inaccessible control sequence so that your^^J
definition will be completed without mixing me up too badly.^^J
If you typed \string\define{\string\cs} the inaccessible control sequence^^J
was defined to be \string\cs, and the rest of your^^J
definition appears as input.}
\newhelp\defbhelp@{I've ignored your definition, because it might^^J
conflict with other uses that are important to me.}
\def\define{\FN@\define@}
\def\define@{\ifcat\noexpand\next\relax
 \expandafter\define@@\else\errhelp\defahelp@                               
 \err@{\string\define\space must be followed by a control
 sequence}\expandafter\def\expandafter\nextii@\fi}                          
\def\undefined@@@@@@@@@@{}
\def\preloaded@@@@@@@@@@{}
\def\next@@@@@@@@@@{}
\def\define@@#1{\ifx#1\relax\errhelp\defbhelp@                              
 \err@{\string#1\space is already defined}\DN@{\DNii@}\else
 \expandafter\ifx\csname\expandafter\eat@\string                            
 #1@@@@@@@@@@\endcsname\undefined@@@@@@@@@@\errhelp\defbhelp@
 \err@{\string#1\space can't be defined}\DN@{\DNii@}\else
 \expandafter\ifx\csname\expandafter\eat@\string#1\endcsname\relax          
 \global\let#1\undefined\DN@{\def#1}\else\errhelp\defbhelp@
 \err@{\string#1\space is already defined}\DN@{\DNii@}\fi
 \fi\fi\next@}

\def\predefine#1#2{\let#1#2}
\def\undefine#1{\let#1\undefined}
\message{page layout,}
\newdimen\captionwidth@
\captionwidth@\hsize
\advance\captionwidth@-1.5in
\def\pagewidth#1{\hsize#1\relax
 \captionwidth@\hsize\advance\captionwidth@-1.5in}
\def\pageheight#1{\vsize#1\relax}
\def\hcorrection#1{\advance\hoffset#1\relax}
\def\vcorrection#1{\advance\voffset#1\relax}
\message{accents/punctuation,}

\let\graveaccent\`
\let\acuteaccent\'
\let\tildeaccent\~
\let\hataccent\^
\let\underscore\_
\let\B\=
\let\D\.
\let\ic@\/
\def\/{\unskip\ic@}
\def\textfonti{\the\textfont\@ne}
\def\t#1#2{{\edef\next@{\the\font}\textfonti\accent"7F \next@#1#2}}
\def~{\unskip\nobreak\ \ignorespaces}
\def\.{.\spacefactor\@m}
\atdef@;{\leavevmode\null;}
\atdef@:{\leavevmode\null:}
\atdef@?{\leavevmode\null?}
\edef\@{\string @}
\def\relaxnext@{\let\next\relax}
\atdef@-{\relaxnext@\leavevmode
 \DN@{\ifx\next-\DN@-{\FN@\nextii@}\else
  \DN@{\leavevmode\hbox{-}}\fi\next@}%
 \DNii@{\ifx\next-\DN@-{\leavevmode\hbox{---}}\else
  \DN@{\leavevmode\hbox{--}}\fi\next@}%
 \FN@\next@}
\def\srdr@{\kern.16667em}
\def\drsr@{\kern.02778em}
\def\sldl@{\drsr@}
\def\dlsl@{\srdr@}
\atdef@"{\unskip\relaxnext@
 \DN@{\ifx\next\space@\DN@. {\FN@\nextii@}\else
  \DN@.{\FN@\nextii@}\fi\next@.}%
 \DNii@{\ifx\next`\DN@`{\FN@\nextiii@}\else
  \ifx\next\lq\DN@\lq{\FN@\nextiii@}\else
  \DN@####1{\FN@\nextiv@}\fi\fi\next@}%
 \def\nextiii@{\ifx\next`\DN@`{\sldl@``}\else\ifx\next\lq
  \DN@\lq{\sldl@``}\else\DN@{\dlsl@`}\fi\fi\next@}%
 \def\nextiv@{\ifx\next'\DN@'{\srdr@''}\else
  \ifx\next\rq\DN@\rq{\srdr@''}\else\DN@{\drsr@'}\fi\fi\next@}%
 \FN@\next@}

\def\textfontii{\the\textfont\tw@}
\def\lbrace@{\delimiter"4266308 }
\def\rbrace@{\delimiter"5267309 }
\def\{{\RIfM@\lbrace@\else{\textfontii f}\spacefactor\@m\fi}
\def\}{\RIfM@\rbrace@\else
 \let\@sf\empty\ifhmode\edef\@sf{\spacefactor\the\spacefactor}\fi
 {\textfontii g}\@sf\relax\fi}
\let\lbrace\{
\let\rbrace\}
\def\AmSTeX{{\textfontii A\kern-.1667em%
  \lower.5ex\hbox{M}\kern-.125emS}-\TeX}
\message{line and page breaks,}
\def\vmodeerr@#1{\Err@{\string#1\space not allowed between paragraphs}}
\def\mathmodeerr@#1{\Err@{\string#1\space not allowed in math mode}}
\def\linebreak{\RIfM@\mathmodeerr@\linebreak\else
 \ifhmode\unskip\unkern\break\else\vmodeerr@\linebreak\fi\fi}

\newskip\saveskip@
\def\allowlinebreak{\RIfM@\mathmodeerr@\allowlinebreak\else
 \ifhmode\saveskip@\lastskip\unskip
 \allowbreak\ifdim\saveskip@>\z@\hskip\saveskip@\fi
 \else\vmodeerr@\allowlinebreak\fi\fi}
\def\nolinebreak{\RIfM@\mathmodeerr@\nolinebreak\else
 \ifhmode\saveskip@\lastskip\unskip
 \nobreak\ifdim\saveskip@>\z@\hskip\saveskip@\fi
 \else\vmodeerr@\nolinebreak\fi\fi}
\def\newline{\relaxnext@
 \DN@{\RIfM@\expandafter\mathmodeerr@\expandafter\newline\else
  \ifhmode\ifx\next\par\else
  \expandafter\unskip\expandafter\null\expandafter\hfill\expandafter\break\fi
  \else
  \expandafter\vmodeerr@\expandafter\newline\fi\fi}%
 \FN@\next@}
\def\dmatherr@#1{\Err@{\string#1\space not allowed in display math mode}}
\def\nondmatherr@#1{\Err@{\string#1\space not allowed in non-display math
 mode}}
\def\onlydmatherr@#1{\Err@{\string#1\space allowed only in display math mode}}
\def\nonmatherr@#1{\Err@{\string#1\space allowed only in math mode}}
\def\mathbreak{\RIfMIfI@\break\else
 \dmatherr@\mathbreak\fi\else\nonmatherr@\mathbreak\fi}
\def\nomathbreak{\RIfMIfI@\nobreak\else
 \dmatherr@\nomathbreak\fi\else\nonmatherr@\nomathbreak\fi}
\def\allowmathbreak{\RIfMIfI@\allowbreak\else
 \dmatherr@\allowmathbreak\fi\else\nonmatherr@\allowmathbreak\fi}
\def\pagebreak{\RIfM@
 \ifinner\nondmatherr@\pagebreak\else\postdisplaypenalty-\@M\fi
 \else\ifvmode\removelastskip\break\else\vadjust{\break}\fi\fi}
\def\nopagebreak{\RIfM@
 \ifinner\nondmatherr@\nopagebreak\else\postdisplaypenalty\@M\fi
 \else\ifvmode\nobreak\else\vadjust{\nobreak}\fi\fi}
\def\nonvmodeerr@#1{\Err@{\string#1\space not allowed within a paragraph
 or in math}}
\def\vnonvmode@#1#2{\relaxnext@\DNii@{\ifx\next\par\DN@{#1}\else
 \DN@{#2}\fi\next@}%
 \ifvmode\DN@{#1}\else
 \DN@{\FN@\nextii@}\fi\next@}
\def\newpage{\vnonvmode@{\vfill\break}{\nonvmodeerr@\newpage}}
\def\smallpagebreak{\vnonvmode@\smallbreak{\nonvmodeerr@\smallpagebreak}}
\def\medpagebreak{\vnonvmode@\medbreak{\nonvmodeerr@\medpagebreak}}
\def\bigpagebreak{\vnonvmode@\bigbreak{\nonvmodeerr@\bigpagebreak}}
\def\NoBlackBoxes{\global\overfullrule\z@}
\def\BlackBoxes{\global\overfullrule5\p@}
\def\Invalid@#1{\def#1{\Err@{\Invalid@@\string#1}}}
\def\Invalid@@{Invalid use of }
\message{figures,}
\Invalid@\caption
\Invalid@\captionwidth
\newdimen\smallcaptionwidth@
\def\topspace{\mid@false\ins@}
\def\midspace{\mid@true\ins@}
\newif\ifmid@
\def\captionfont@{}
\def\ins@#1{\relaxnext@\allowbreak
 \smallcaptionwidth@\captionwidth@\gdef\thespace@{#1}%
 \DN@{\ifx\next\space@\DN@. {\FN@\nextii@}\else
  \DN@.{\FN@\nextii@}\fi\next@.}%
 \DNii@{\ifx\next\caption\DN@\caption{\FN@\nextiii@}%
  \else\let\next@\nextiv@\fi\next@}%
 \def\nextiv@{\vnonvmode@
  {\ifmid@\expandafter\midinsert\else\expandafter\topinsert\fi
   \vbox to\thespace@{}\endinsert}
  {\ifmid@\nonvmodeerr@\midspace\else\nonvmodeerr@\topspace\fi}}%
 \def\nextiii@{\ifx\next\captionwidth\expandafter\nextv@
  \else\expandafter\nextvi@\fi}%
 \def\nextv@\captionwidth##1##2{\smallcaptionwidth@##1\relax\nextvi@{##2}}%
 \def\nextvi@##1{\def\thecaption@{\captionfont@##1}%
  \DN@{\ifx\next\space@\DN@. {\FN@\nextvii@}\else
   \DN@.{\FN@\nextvii@}\fi\next@.}%
  \FN@\next@}%
 \def\nextvii@{\vnonvmode@
  {\ifmid@\expandafter\midinsert\else
  \expandafter\topinsert\fi\vbox to\thespace@{}\nobreak\smallskip
  \setboxz@h{\noindent\ignorespaces\thecaption@\unskip}%
  \ifdim\wdz@>\smallcaptionwidth@\centerline{\vbox{\hsize\smallcaptionwidth@
   \noindent\ignorespaces\thecaption@\unskip}}%
  \else\centerline{\boxz@}\fi\endinsert}
  {\ifmid@\nonvmodeerr@\midspace
  \else\nonvmodeerr@\topspace\fi}}%
 \FN@\next@}
\message{comments,}
\def\newcodes@{\catcode`\\12\catcode`\{12\catcode`\}12\catcode`\#12%
 \catcode`\%12\relax}
\def\oldcodes@{\catcode`\\0\catcode`\{1\catcode`\}2\catcode`\#6%
 \catcode`\%14\relax}
\def\comment{\newcodes@\endlinechar=10 \comment@}
{\lccode`\0=`\\
\lowercase{\gdef\comment@#1^^J{\comment@@#10endcomment\comment@@@}%
\gdef\comment@@#10endcomment{\FN@\comment@@@}%
\gdef\comment@@@#1\comment@@@{\ifx\next\comment@@@\let\next\comment@
 \else\def\next{\oldcodes@\endlinechar=`\^^M\relax}%
 \fi\next}}}
\def\pr@m@s{\ifx'\next\DN@##1{\prim@s}\else\let\next@\egroup\fi\next@}
\def\prime{{\null\prime@\null}}
\mathchardef\prime@="0230
\let\dsize\displaystyle

\let\ssize\scriptstyle

\message{math spacing,}
\def\,{\RIfM@\mskip\thinmuskip\relax\else\kern.16667em\fi}
\def\!{\RIfM@\mskip-\thinmuskip\relax\else\kern-.16667em\fi}
\let\thinspace\,
\let\negthinspace\!
\def\medspace{\RIfM@\mskip\medmuskip\relax\else\kern.222222em\fi}
\def\negmedspace{\RIfM@\mskip-\medmuskip\relax\else\kern-.222222em\fi}
\def\thickspace{\RIfM@\mskip\thickmuskip\relax\else\kern.27777em\fi}
\let\;\thickspace
\def\negthickspace{\RIfM@\mskip-\thickmuskip\relax\else
 \kern-.27777em\fi}
\atdef@,{\RIfM@\mskip.1\thinmuskip\else\leavevmode\null,\fi}
\atdef@!{\RIfM@\mskip-.1\thinmuskip\else\leavevmode\null!\fi}
\atdef@.{\RIfM@&&\else\leavevmode.\spacefactor3000 \fi}
\def\and{\DOTSB\;\mathchar"3026 \;}

\message{fractions,}
\def\frac#1#2{{#1\over#2}}

\newdimen\ex@
\ex@.2326ex
\Invalid@\thickness
\def\thickfrac{\relaxnext@
 \DN@{\ifx\next\thickness\let\next@\nextii@\else
 \DN@{\nextii@\thickness1}\fi\next@}%
 \DNii@\thickness##1##2##3{{##2\above##1\ex@##3}}%
 \FN@\next@}

\def\thickfracwithdelims#1#2{\relaxnext@\def\ldelim@{#1}\def\rdelim@{#2}%
 \DN@{\ifx\next\thickness\let\next@\nextii@\else
 \DN@{\nextii@\thickness1}\fi\next@}%
 \DNii@\thickness##1##2##3{{##2\abovewithdelims
 \ldelim@\rdelim@##1\ex@##3}}%
 \FN@\next@}

\def\:{\nobreak\hskip.1111em\mathpunct{}\nonscript\mkern-\thinmuskip{:}\hskip
 .3333emplus.0555em\relax}
\def\snug{\unskip\kern-\mathsurround}
\message{smash commands,}
\def\topsmash{\top@true\bot@false\smash@}
\def\botsmash{\top@false\bot@true\smash@}
\newif\iftop@
\newif\ifbot@
\def\smash{\top@true\bot@true\smash@}
\def\smash@{\RIfM@\expandafter\mathpalette\expandafter\mathsm@sh\else
 \expandafter\makesm@sh\fi}
\def\finsm@sh{\iftop@\ht\z@\z@\fi\ifbot@\dp\z@\z@\fi\leavevmode\boxz@}
\message{large operator symbols,}
\def\LimitsOnSums{\global\let\slimits@\displaylimits}
\def\NoLimitsOnSums{\global\let\slimits@\nolimits}
\LimitsOnSums
\mathchardef\coprod@="1360       \def\coprod{\DOTSB\coprod@\slimits@}
\mathchardef\bigvee@="1357       \def\bigvee{\DOTSB\bigvee@\slimits@}
\mathchardef\bigwedge@="1356     \def\bigwedge{\DOTSB\bigwedge@\slimits@}
\mathchardef\biguplus@="1355     \def\biguplus{\DOTSB\biguplus@\slimits@}
\mathchardef\bigcap@="1354       \def\bigcap{\DOTSB\bigcap@\slimits@}
\mathchardef\bigcup@="1353       \def\bigcup{\DOTSB\bigcup@\slimits@}
\mathchardef\prod@="1351         \def\prod{\DOTSB\prod@\slimits@}
\mathchardef\sum@="1350          \def\sum{\DOTSB\sum@\slimits@}
\mathchardef\bigotimes@="134E    \def\bigotimes{\DOTSB\bigotimes@\slimits@}
\mathchardef\bigoplus@="134C     \def\bigoplus{\DOTSB\bigoplus@\slimits@}
\mathchardef\bigodot@="134A      \def\bigodot{\DOTSB\bigodot@\slimits@}
\mathchardef\bigsqcup@="1346     \def\bigsqcup{\DOTSB\bigsqcup@\slimits@}
\message{integrals,}
\def\LimitsOnInts{\global\let\ilimits@\displaylimits}
\def\NoLimitsOnInts{\global\let\ilimits@\nolimits}
\NoLimitsOnInts
\def\int{\DOTSI\intop\ilimits@}
\def\oint{\DOTSI\ointop\ilimits@}
\def\intic@{\mathchoice{\hskip.5em}{\hskip.4em}{\hskip.4em}{\hskip.4em}}
\def\negintic@{\mathchoice
 {\hskip-.5em}{\hskip-.4em}{\hskip-.4em}{\hskip-.4em}}
\def\intkern@{\mathchoice{\!\!\!}{\!\!}{\!\!}{\!\!}}
\def\intdots@{\mathchoice{\plaincdots@}
 {{\cdotp}\mkern1.5mu{\cdotp}\mkern1.5mu{\cdotp}}
 {{\cdotp}\mkern1mu{\cdotp}\mkern1mu{\cdotp}}
 {{\cdotp}\mkern1mu{\cdotp}\mkern1mu{\cdotp}}}
\newcount\intno@
\def\iint{\DOTSI\intno@\tw@\FN@\ints@}
\def\iiint{\DOTSI\intno@\thr@@\FN@\ints@}
\def\iiiint{\DOTSI\intno@4 \FN@\ints@}
\def\idotsint{\DOTSI\intno@\z@\FN@\ints@}
\def\ints@{\findlimits@\ints@@}
\newif\iflimtoken@
\newif\iflimits@
\def\findlimits@{\limtoken@true\ifx\next\limits\limits@true
 \else\ifx\next\nolimits\limits@false\else
 \limtoken@false\ifx\ilimits@\nolimits\limits@false\else
 \ifinner\limits@false\else\limits@true\fi\fi\fi\fi}
\def\multint@{\int\ifnum\intno@=\z@\intdots@                                
 \else\intkern@\fi                                                          
 \ifnum\intno@>\tw@\int\intkern@\fi                                         
 \ifnum\intno@>\thr@@\int\intkern@\fi                                       
 \int}                                                                      
\def\multintlimits@{\intop\ifnum\intno@=\z@\intdots@\else\intkern@\fi
 \ifnum\intno@>\tw@\intop\intkern@\fi
 \ifnum\intno@>\thr@@\intop\intkern@\fi\intop}
\def\ints@@{\iflimtoken@                                                    
 \def\ints@@@{\iflimits@\negintic@\mathop{\intic@\multintlimits@}\limits    
  \else\multint@\nolimits\fi                                                
  \eat@}                                                                    
 \else                                                                      
 \def\ints@@@{\iflimits@\negintic@
  \mathop{\intic@\multintlimits@}\limits\else
  \multint@\nolimits\fi}\fi\ints@@@}
\def\LimitsOnNames{\global\let\nlimits@\displaylimits}
\def\NoLimitsOnNames{\global\let\nlimits@\nolimits@}
\LimitsOnNames
\def\nolimits@{\relaxnext@
 \DN@{\ifx\next\limits\DN@\limits{\nolimits}\else
  \let\next@\nolimits\fi\next@}%
 \FN@\next@}
\message{operator names,}
\def\newmcodes@{\mathcode`\'"27\mathcode`\*"2A\mathcode`\."613A%
 \mathcode`\-"2D\mathcode`\/"2F\mathcode`\:"603A }
\def\operatorname#1{\mathop{\newmcodes@\kern\z@\fam\z@#1}\nolimits@}
\def\operatornamewithlimits#1{\mathop{\newmcodes@\kern\z@\fam\z@#1}\nlimits@}
\def\qopname@#1{\mathop{\fam\z@#1}\nolimits@}
\def\qopnamewl@#1{\mathop{\fam\z@#1}\nlimits@}
\def\arccos{\qopname@{arccos}}
\def\arcsin{\qopname@{arcsin}}
\def\arctan{\qopname@{arctan}}
\def\arg{\qopname@{arg}}
\def\cos{\qopname@{cos}}
\def\cosh{\qopname@{cosh}}
\def\cot{\qopname@{cot}}
\def\coth{\qopname@{coth}}
\def\csc{\qopname@{csc}}
\def\deg{\qopname@{deg}}
\def\det{\qopnamewl@{det}}
\def\dim{\qopname@{dim}}
\def\exp{\qopname@{exp}}
\def\gcd{\qopnamewl@{gcd}}
\def\hom{\qopname@{hom}}
\def\inf{\qopnamewl@{inf}}
\def\injlim{\qopnamewl@{inj\,lim}}
\def\ker{\qopname@{ker}}
\def\lg{\qopname@{lg}}
\def\lim{\qopnamewl@{lim}}
\def\liminf{\qopnamewl@{lim\,inf}}
\def\limsup{\qopnamewl@{lim\,sup}}
\def\ln{\qopname@{ln}}
\def\log{\qopname@{log}}
\def\max{\qopnamewl@{max}}
\def\min{\qopnamewl@{min}}
\def\Pr{\qopnamewl@{Pr}}
\def\projlim{\qopnamewl@{proj\,lim}}
\def\sec{\qopname@{sec}}
\def\sin{\qopname@{sin}}
\def\sinh{\qopname@{sinh}}
\def\sup{\qopnamewl@{sup}}
\def\tan{\qopname@{tan}}
\def\tanh{\qopname@{tanh}}
\def\varinjlim{\mathop{\vtop{\ialign{##\crcr
 \hfil\rm lim\hfil\crcr\noalign{\nointerlineskip}\rightarrowfill\crcr
 \noalign{\nointerlineskip\kern-\ex@}\crcr}}}}
\def\varprojlim{\mathop{\vtop{\ialign{##\crcr
 \hfil\rm lim\hfil\crcr\noalign{\nointerlineskip}\leftarrowfill\crcr
 \noalign{\nointerlineskip\kern-\ex@}\crcr}}}}
\def\varliminf{\mathop{\underline{\vrule height\z@ depth.2exwidth\z@
 \hbox{\rm lim}}}}

\newdimen\buffer@
\buffer@\fontdimen13 \tenex
\newdimen\buffer
\buffer\buffer@

\def\ResetBuffer{\fontdimen13 \tenex\buffer@\global\buffer\buffer@}
\def\shave#1{\mathop{\hbox{$\m@th\fontdimen13 \tenex\z@                     
 \displaystyle{#1}$}}\fontdimen13 \tenex\buffer}

\message{multilevel sub/superscripts,}
\Invalid@\\
\def\Let@{\relax\iffalse{\fi\let\\=\cr\iffalse}\fi}
\Invalid@\vspace
\def\vspace@{\def\vspace##1{\crcr\noalign{\vskip##1\relax}}}
\def\multilimits@{\bgroup\vspace@\Let@
 \baselineskip\fontdimen10 \scriptfont\tw@
 \advance\baselineskip\fontdimen12 \scriptfont\tw@
 \lineskip\thr@@\fontdimen8 \scriptfont\thr@@
 \lineskiplimit\lineskip
 \vbox\bgroup\ialign\bgroup\hfil$\m@th\scriptstyle{##}$\hfil\crcr}
\def\Sb{_\multilimits@}
\def\endSb{\crcr\egroup\egroup\egroup}
\def\Sp{^\multilimits@}

\def\spreadlines#1{\RIfMIfI@\onlydmatherr@\spreadlines\else
 \openup#1\relax\fi\else\onlydmatherr@\spreadlines\fi}
\def\Mathstrut@{\copy\Mathstrutbox@}
\newbox\Mathstrutbox@
\setbox\Mathstrutbox@\null
\setboxz@h{$\m@th($}
\ht\Mathstrutbox@\ht\z@
\dp\Mathstrutbox@\dp\z@
\message{matrices,}
\newdimen\spreadmlines@
\def\spreadmatrixlines#1{\RIfMIfI@
 \onlydmatherr@\spreadmatrixlines\else
 \spreadmlines@#1\relax\fi\else\onlydmatherr@\spreadmatrixlines\fi}
\def\matrix{\null\,\vcenter\bgroup\Let@\vspace@
 \normalbaselines\openup\spreadmlines@\ialign
 \bgroup\hfil$\m@th##$\hfil&&\quad\hfil$\m@th##$\hfil\crcr
 \Mathstrut@\crcr\noalign{\kern-\baselineskip}}
\def\endmatrix{\crcr\Mathstrut@\crcr\noalign{\kern-\baselineskip}\egroup
 \egroup\,}
\def\format{\crcr\egroup\iffalse{\fi\ifnum`}=0 \fi\format@}
\newtoks\hashtoks@
\hashtoks@{#}
\def\format@#1\\{\def\preamble@{#1}%
 \def\l{$\m@th\the\hashtoks@$\hfil}%
 \def\c{\hfil$\m@th\the\hashtoks@$\hfil}%
 \def\r{\hfil$\m@th\the\hashtoks@$}%
 \edef\preamble@@{\preamble@}\ifnum`{=0 \fi\iffalse}\fi
 \ialign\bgroup\span\preamble@@\crcr}
\def\smallmatrix{\null\,\vcenter\bgroup\vspace@\Let@
 \baselineskip9\ex@\lineskip\ex@
 \ialign\bgroup\hfil$\m@th\scriptstyle{##}$\hfil&&\thickspace\hfil
 $\m@th\scriptstyle{##}$\hfil\crcr}
\def\endsmallmatrix{\crcr\egroup\egroup\,}

\newmuskip\dotsspace@
\dotsspace@1.5mu
\def\strip@#1 {#1}
\def\spacehdots#1\for#2{\multispan{#2}\xleaders
 \hbox{$\m@th\mkern\strip@#1 \dotsspace@.\mkern\strip@#1 \dotsspace@$}\hfill}
\def\hdotsfor#1{\spacehdots\@ne\for{#1}}
\def\multispan@#1{\omit\mscount#1\unskip\loop\ifnum\mscount>\@ne\sp@n\repeat}
\def\spaceinnerhdots#1\for#2\after#3{\multispan@{\strip@#2 }#3\xleaders
 \hbox{$\m@th\mkern\strip@#1 \dotsspace@.\mkern\strip@#1 \dotsspace@$}\hfill}
\def\innerhdotsfor#1\after#2{\spaceinnerhdots\@ne\for#1\after{#2}}
\def\cases{\bgroup\spreadmlines@\jot\left\{\,\matrix\format\l&\quad\l\\}
\def\endcases{\endmatrix\right.\egroup}
\message{multiline displays,}
\newif\ifinany@
\newif\ifinalign@
\newif\ifingather@
\def\strut@{\copy\strutbox@}
\newbox\strutbox@
\setbox\strutbox@\hbox{\vrule height8\p@ depth3\p@ width\z@}
\def\topaligned{\null\,\vtop\aligned@}
\def\botaligned{\null\,\vbox\aligned@}
\def\aligned{\null\,\vcenter\aligned@}
\def\aligned@{\bgroup\vspace@\Let@
 \ifinany@\else\openup\jot\fi\ialign
 \bgroup\hfil\strut@$\m@th\displaystyle{##}$&
 $\m@th\displaystyle{{}##}$\hfil\crcr}
\def\endaligned{\crcr\egroup\egroup}

\def\alignedat#1{\null\,\vcenter\bgroup\doat@{#1}\vspace@\Let@
 \ifinany@\else\openup\jot\fi\ialign\bgroup\span\preamble@@\crcr}
\newcount\atcount@
\def\doat@#1{\toks@{\hfil\strut@$\m@th
 \displaystyle{\the\hashtoks@}$&$\m@th\displaystyle
 {{}\the\hashtoks@}$\hfil}
 \atcount@#1\relax\advance\atcount@\m@ne                                    
 \loop\ifnum\atcount@>\z@\toks@=\expandafter{\the\toks@&\hfil$\m@th
 \displaystyle{\the\hashtoks@}$&$\m@th
 \displaystyle{{}\the\hashtoks@}$\hfil}\advance
  \atcount@\m@ne\repeat                                                     
 \xdef\preamble@{\the\toks@}\xdef\preamble@@{\preamble@}}

\def\gathered{\null\,\vcenter\bgroup\vspace@\Let@
 \ifinany@\else\openup\jot\fi\ialign
 \bgroup\hfil\strut@$\m@th\displaystyle{##}$\hfil\crcr}
\def\endgathered{\crcr\egroup\egroup}
\newif\iftagsleft@
\def\TagsOnLeft{\global\tagsleft@true}
\def\TagsOnRight{\global\tagsleft@false}
\TagsOnLeft
\newif\ifmathtags@
\def\TagsAsMath{\global\mathtags@true}
\def\TagsAsText{\global\mathtags@false}
\TagsAsText
\def\tagform@#1{\hbox{\rm(\ignorespaces#1\unskip)}}
\def\thetag{\leavevmode\tagform@}
\def\tag#1$${\iftagsleft@\leqno\else\eqno\fi                                
 \maketag@#1\maketag@                                                       
 $$}                                                                        
\def\maketag@{\FN@\maketag@@}
\def\maketag@@{\ifx\next"\expandafter\maketag@@@\else\expandafter\maketag@@@@
 \fi}
\def\maketag@@@"#1"#2\maketag@{\hbox{\rm#1}}                                
\def\maketag@@@@#1\maketag@{\ifmathtags@\tagform@{$\m@th#1$}\else
 \tagform@{#1}\fi}
\interdisplaylinepenalty\@M
\def\allowdisplaybreaks{\RIfMIfI@
 \onlydmatherr@\allowdisplaybreaks\else
 \interdisplaylinepenalty\z@\fi\else\onlydmatherr@\allowdisplaybreaks\fi}
\Invalid@\allowdisplaybreak
\Invalid@\displaybreak
\Invalid@\intertext
\def\allowdisplaybreak@{\def\allowdisplaybreak{\crcr\noalign{\allowbreak}}}
\def\displaybreak@{\def\displaybreak{\crcr\noalign{\break}}}
\def\intertext@{\def\intertext##1{\crcr\noalign{%
 \penalty\postdisplaypenalty \vskip\belowdisplayskip
 \vbox{\normalbaselines\noindent##1}%
 \penalty\predisplaypenalty \vskip\abovedisplayskip}}}
\newskip\centering@
\centering@\z@ plus\@m\p@
\def\align{\relax\ifingather@\DN@{\csname align (in
  \string\gather)\endcsname}\else
 \ifmmode\ifinner\DN@{\onlydmatherr@\align}\else
  \let\next@\align@\fi
 \else\DN@{\onlydmatherr@\align}\fi\fi\next@}
\newhelp\andhelp@
{An extra & here is so disastrous that you should probably exit^^J
and fix things up.}
\newif\iftag@
\newcount\and@
\def\align@{\inalign@true\inany@true
 \vspace@\allowdisplaybreak@\displaybreak@\intertext@
 \def\tag{\global\tag@true\ifnum\and@=\z@\DN@{&&}\else
          \DN@{&}\fi\next@}%
 \iftagsleft@\DN@{\csname align \endcsname}\else
  \DN@{\csname align \space\endcsname}\fi\next@}
\def\Tag@{\iftag@\else\errhelp\andhelp@\err@{Extra & on this line}\fi}
\newdimen\lwidth@
\newdimen\rwidth@
\newdimen\maxlwidth@
\newdimen\maxrwidth@
\newdimen\totwidth@
\def\measure@#1\endalign{\lwidth@\z@\rwidth@\z@\maxlwidth@\z@\maxrwidth@\z@
 \global\and@\z@                                                            
 \setbox@ne\vbox                                                            
  {\everycr{\noalign{\global\tag@false\global\and@\z@}}\Let@                
  \halign{\setboxz@h{$\m@th\displaystyle{\@lign##}$}
   \global\lwidth@\wdz@                                                     
   \ifdim\lwidth@>\maxlwidth@\global\maxlwidth@\lwidth@\fi                  
   \global\advance\and@\@ne                                                 
   &\setboxz@h{$\m@th\displaystyle{{}\@lign##}$}\global\rwidth@\wdz@        
   \ifdim\rwidth@>\maxrwidth@\global\maxrwidth@\rwidth@\fi                  
   \global\advance\and@\@ne                                                
   &\Tag@
   \eat@{##}\crcr#1\crcr}}
 \totwidth@\maxlwidth@\advance\totwidth@\maxrwidth@}                       
\def\displ@y@{\global\dt@ptrue\openup\jot
 \everycr{\noalign{\global\tag@false\global\and@\z@\ifdt@p\global\dt@pfalse
 \vskip-\lineskiplimit\vskip\normallineskiplimit\else
 \penalty\interdisplaylinepenalty\fi}}}
\def\black@#1{\noalign{\ifdim#1>\displaywidth
 \dimen@\prevdepth\nointerlineskip                                          
 \vskip-\ht\strutbox@\vskip-\dp\strutbox@                                   
 \vbox{\noindent\hbox to#1{\strut@\hfill}}
 \prevdepth\dimen@                                                          
 \fi}}
\expandafter\def\csname align \space\endcsname#1\endalign
 {\measure@#1\endalign\global\and@\z@                                       
 \ifingather@\everycr{\noalign{\global\and@\z@}}\else\displ@y@\fi           
 \Let@\tabskip\centering@                                                   
 \halign to\displaywidth
  {\hfil\strut@\setboxz@h{$\m@th\displaystyle{\@lign##}$}
  \global\lwidth@\wdz@\boxz@\global\advance\and@\@ne                        
  \tabskip\z@skip                                                           
  &\setboxz@h{$\m@th\displaystyle{{}\@lign##}$}
  \global\rwidth@\wdz@\boxz@\hfill\global\advance\and@\@ne                  
  \tabskip\centering@                                                       
  &\setboxz@h{\@lign\strut@\maketag@##\maketag@}
  \dimen@\displaywidth\advance\dimen@-\totwidth@
  \divide\dimen@\tw@\advance\dimen@\maxrwidth@\advance\dimen@-\rwidth@     
  \ifdim\dimen@<\tw@\wdz@\llap{\vtop{\normalbaselines\null\boxz@}}
  \else\llap{\boxz@}\fi                                                    
  \tabskip\z@skip                                                          
  \crcr#1\crcr                                                             
  \black@\totwidth@}}                                                      
\newdimen\lineht@
\expandafter\def\csname align \endcsname#1\endalign{\measure@#1\endalign
 \global\and@\z@
 \ifdim\totwidth@>\displaywidth\let\displaywidth@\totwidth@\else
  \let\displaywidth@\displaywidth\fi                                        
 \ifingather@\everycr{\noalign{\global\and@\z@}}\else\displ@y@\fi
 \Let@\tabskip\centering@\halign to\displaywidth
  {\hfil\strut@\setboxz@h{$\m@th\displaystyle{\@lign##}$}%
  \global\lwidth@\wdz@\global\lineht@\ht\z@                                 
  \boxz@\global\advance\and@\@ne
  \tabskip\z@skip&\setboxz@h{$\m@th\displaystyle{{}\@lign##}$}%
  \global\rwidth@\wdz@\ifdim\ht\z@>\lineht@\global\lineht@\ht\z@\fi         
  \boxz@\hfil\global\advance\and@\@ne
  \tabskip\centering@&\kern-\displaywidth@                                  
  \setboxz@h{\@lign\strut@\maketag@##\maketag@}%
  \dimen@\displaywidth\advance\dimen@-\totwidth@
  \divide\dimen@\tw@\advance\dimen@\maxlwidth@\advance\dimen@-\lwidth@
  \ifdim\dimen@<\tw@\wdz@
   \rlap{\vbox{\normalbaselines\boxz@\vbox to\lineht@{}}}\else
   \rlap{\boxz@}\fi
  \tabskip\displaywidth@\crcr#1\crcr\black@\totwidth@}}
\expandafter\def\csname align (in \string\gather)\endcsname
  #1\endalign{\vcenter{\align@#1\endalign}}
\Invalid@\endalign
\newif\ifxat@
\def\alignat{\RIfMIfI@\DN@{\onlydmatherr@\alignat}\else
 \DN@{\csname alignat \endcsname}\fi\else
 \DN@{\onlydmatherr@\alignat}\fi\next@}
\newif\ifmeasuring@
\newbox\savealignat@
\expandafter\def\csname alignat \endcsname#1#2\endalignat                   
 {\inany@true\xat@false
 \def\tag{\global\tag@true\count@#1\relax\multiply\count@\tw@
  \xdef\tag@{}\loop\ifnum\count@>\and@\xdef\tag@{&\tag@}\advance\count@\m@ne
  \repeat\tag@}%
 \vspace@\allowdisplaybreak@\displaybreak@\intertext@
 \displ@y@\measuring@true                                                   
 \setbox\savealignat@\hbox{$\m@th\displaystyle\Let@
  \attag@{#1}
  \vbox{\halign{\span\preamble@@\crcr#2\crcr}}$}%
 \measuring@false                                                           
 \Let@\attag@{#1}
 \tabskip\centering@\halign to\displaywidth
  {\span\preamble@@\crcr#2\crcr                                             
  \black@{\wd\savealignat@}}}                                               
\Invalid@\endalignat
\def\xalignat{\RIfMIfI@
 \DN@{\onlydmatherr@\xalignat}\else
 \DN@{\csname xalignat \endcsname}\fi\else
 \DN@{\onlydmatherr@\xalignat}\fi\next@}
\expandafter\def\csname xalignat \endcsname#1#2\endxalignat
 {\inany@true\xat@true
 \def\tag{\global\tag@true\def\tag@{}\count@#1\relax\multiply\count@\tw@
  \loop\ifnum\count@>\and@\xdef\tag@{&\tag@}\advance\count@\m@ne\repeat\tag@}%
 \vspace@\allowdisplaybreak@\displaybreak@\intertext@
 \displ@y@\measuring@true\setbox\savealignat@\hbox{$\m@th\displaystyle\Let@
 \attag@{#1}\vbox{\halign{\span\preamble@@\crcr#2\crcr}}$}%
 \measuring@false\Let@
 \attag@{#1}\tabskip\centering@\halign to\displaywidth
 {\span\preamble@@\crcr#2\crcr\black@{\wd\savealignat@}}}
\def\attag@#1{\let\Maketag@\maketag@\let\TAG@\Tag@                          
 \let\Tag@=0\let\maketag@=0
 \ifmeasuring@\def\llap@##1{\setboxz@h{##1}\hbox to\tw@\wdz@{}}%
  \def\rlap@##1{\setboxz@h{##1}\hbox to\tw@\wdz@{}}\else
  \let\llap@\llap\let\rlap@\rlap\fi                                         
 \toks@{\hfil\strut@$\m@th\displaystyle{\@lign\the\hashtoks@}$\tabskip\z@skip
  \global\advance\and@\@ne&$\m@th\displaystyle{{}\@lign\the\hashtoks@}$\hfil
  \ifxat@\tabskip\centering@\fi\global\advance\and@\@ne}
 \iftagsleft@
  \toks@@{\tabskip\centering@&\Tag@\kern-\displaywidth
   \rlap@{\@lign\maketag@\the\hashtoks@\maketag@}%
   \global\advance\and@\@ne\tabskip\displaywidth}\else
  \toks@@{\tabskip\centering@&\Tag@\llap@{\@lign\maketag@
   \the\hashtoks@\maketag@}\global\advance\and@\@ne\tabskip\z@skip}\fi      
 \atcount@#1\relax\advance\atcount@\m@ne
 \loop\ifnum\atcount@>\z@
 \toks@=\expandafter{\the\toks@&\hfil$\m@th\displaystyle{\@lign
  \the\hashtoks@}$\global\advance\and@\@ne
  \tabskip\z@skip&$\m@th\displaystyle{{}\@lign\the\hashtoks@}$\hfil\ifxat@
  \tabskip\centering@\fi\global\advance\and@\@ne}\advance\atcount@\m@ne
 \repeat                                                                    
 \xdef\preamble@{\the\toks@\the\toks@@}
 \xdef\preamble@@{\preamble@}
 \let\maketag@\Maketag@\let\Tag@\TAG@}                                      
\Invalid@\endxalignat
\def\xxalignat{\RIfMIfI@
 \DN@{\onlydmatherr@\xxalignat}\else\DN@{\csname xxalignat
  \endcsname}\fi\else
 \DN@{\onlydmatherr@\xxalignat}\fi\next@}
\expandafter\def\csname xxalignat \endcsname#1#2\endxxalignat{\inany@true
 \vspace@\allowdisplaybreak@\displaybreak@\intertext@
 \displ@y\setbox\savealignat@\hbox{$\m@th\displaystyle\Let@
 \xxattag@{#1}\vbox{\halign{\span\preamble@@\crcr#2\crcr}}$}%
 \Let@\xxattag@{#1}\tabskip\z@skip\halign to\displaywidth
 {\span\preamble@@\crcr#2\crcr\black@{\wd\savealignat@}}}
\def\xxattag@#1{\toks@{\tabskip\z@skip\hfil\strut@
 $\m@th\displaystyle{\the\hashtoks@}$&%
 $\m@th\displaystyle{{}\the\hashtoks@}$\hfil\tabskip\centering@&}%
 \atcount@#1\relax\advance\atcount@\m@ne\loop\ifnum\atcount@>\z@
 \toks@=\expandafter{\the\toks@&\hfil$\m@th\displaystyle{\the\hashtoks@}$%
  \tabskip\z@skip&$\m@th\displaystyle{{}\the\hashtoks@}$\hfil
  \tabskip\centering@}\advance\atcount@\m@ne\repeat
 \xdef\preamble@{\the\toks@\tabskip\z@skip}\xdef\preamble@@{\preamble@}}
\Invalid@\endxxalignat
\newdimen\gwidth@
\newdimen\gmaxwidth@
\def\gmeasure@#1\endgather{\gwidth@\z@\gmaxwidth@\z@\setbox@ne\vbox{\Let@
 \halign{\setboxz@h{$\m@th\displaystyle{##}$}\global\gwidth@\wdz@
 \ifdim\gwidth@>\gmaxwidth@\global\gmaxwidth@\gwidth@\fi
 &\eat@{##}\crcr#1\crcr}}}
\def\gather{\RIfMIfI@\DN@{\onlydmatherr@\gather}\else
 \ingather@true\inany@true\def\tag{&}%
 \vspace@\allowdisplaybreak@\displaybreak@\intertext@
 \displ@y\Let@
 \iftagsleft@\DN@{\csname gather \endcsname}\else
  \DN@{\csname gather \space\endcsname}\fi\fi
 \else\DN@{\onlydmatherr@\gather}\fi\next@}
\expandafter\def\csname gather \space\endcsname#1\endgather
 {\gmeasure@#1\endgather\tabskip\centering@
 \halign to\displaywidth{\hfil\strut@\setboxz@h{$\m@th\displaystyle{##}$}%
 \global\gwidth@\wdz@\boxz@\hfil&
 \setboxz@h{\strut@{\maketag@##\maketag@}}%
 \dimen@\displaywidth\advance\dimen@-\gwidth@
 \ifdim\dimen@>\tw@\wdz@\llap{\boxz@}\else
 \llap{\vtop{\normalbaselines\null\boxz@}}\fi
 \tabskip\z@skip\crcr#1\crcr\black@\gmaxwidth@}}
\newdimen\glineht@
\expandafter\def\csname gather \endcsname#1\endgather{\gmeasure@#1\endgather
 \ifdim\gmaxwidth@>\displaywidth\let\gdisplaywidth@\gmaxwidth@\else
 \let\gdisplaywidth@\displaywidth\fi\tabskip\centering@\halign to\displaywidth
 {\hfil\strut@\setboxz@h{$\m@th\displaystyle{##}$}%
 \global\gwidth@\wdz@\global\glineht@\ht\z@\boxz@\hfil&\kern-\gdisplaywidth@
 \setboxz@h{\strut@{\maketag@##\maketag@}}%
 \dimen@\displaywidth\advance\dimen@-\gwidth@
 \ifdim\dimen@>\tw@\wdz@\rlap{\boxz@}\else
 \rlap{\vbox{\normalbaselines\boxz@\vbox to\glineht@{}}}\fi
 \tabskip\gdisplaywidth@\crcr#1\crcr\black@\gmaxwidth@}}
\newif\ifctagsplit@
\def\CenteredTagsOnSplits{\global\ctagsplit@true}
\def\TopOrBottomTagsOnSplits{\global\ctagsplit@false}
\TopOrBottomTagsOnSplits
\def\split{\relax\ifinany@\let\next@\insplit@\else
 \ifmmode\ifinner\def\next@{\onlydmatherr@\split}\else
 \let\next@\outsplit@\fi\else
 \def\next@{\onlydmatherr@\split}\fi\fi\next@}
\def\insplit@{\global\setbox\z@\vbox\bgroup\vspace@\Let@\ialign\bgroup
 \hfil\strut@$\m@th\displaystyle{##}$&$\m@th\displaystyle{{}##}$\hfill\crcr}
\def\endsplit{\crcr\egroup\egroup\iftagsleft@\expandafter\lendsplit@\else
 \expandafter\rendsplit@\fi}
\def\rendsplit@{\global\setbox9 \vbox
 {\unvcopy\z@\global\setbox8 \lastbox\unskip}
 \setbox@ne\hbox{\unhcopy8 \unskip\global\setbox\tw@\lastbox
 \unskip\global\setbox\thr@@\lastbox}
 \global\setbox7 \hbox{\unhbox\tw@\unskip}
 \ifinalign@\ifctagsplit@                                                   
  \gdef\split@{\hbox to\wd\thr@@{}&
   \vcenter{\vbox{\moveleft\wd\thr@@\boxz@}}}
 \else\gdef\split@{&\vbox{\moveleft\wd\thr@@\box9}\crcr
  \box\thr@@&\box7}\fi                                                      
 \else                                                                      
  \ifctagsplit@\gdef\split@{\vcenter{\boxz@}}\else
  \gdef\split@{\box9\crcr\hbox{\box\thr@@\box7}}\fi
 \fi
 \split@}                                                                   
\def\lendsplit@{\global\setbox9\vtop{\unvcopy\z@}
 \setbox@ne\vbox{\unvcopy\z@\global\setbox8\lastbox}
 \setbox@ne\hbox{\unhcopy8\unskip\setbox\tw@\lastbox
  \unskip\global\setbox\thr@@\lastbox}
 \ifinalign@\ifctagsplit@                                                   
  \gdef\split@{\hbox to\wd\thr@@{}&
  \vcenter{\vbox{\moveleft\wd\thr@@\box9}}}
  \else                                                                     
  \gdef\split@{\hbox to\wd\thr@@{}&\vbox{\moveleft\wd\thr@@\box9}}\fi
 \else
  \ifctagsplit@\gdef\split@{\vcenter{\box9}}\else
  \gdef\split@{\box9}\fi
 \fi\split@}
\def\outsplit@#1$${\align\insplit@#1\endalign$$}
\newdimen\multlinegap@
\multlinegap@1em
\newdimen\multlinetaggap@
\multlinetaggap@1em
\def\MultlineGap#1{\global\multlinegap@#1\relax}
\def\multlinegap#1{\RIfMIfI@\onlydmatherr@\multlinegap\else
 \multlinegap@#1\relax\fi\else\onlydmatherr@\multlinegap\fi}
\def\nomultlinegap{\multlinegap{\z@}}
\def\multline{\RIfMIfI@
 \DN@{\onlydmatherr@\multline}\else
 \DN@{\multline@}\fi\else
 \DN@{\onlydmatherr@\multline}\fi\next@}
\newif\iftagin@
\def\tagin@#1{\tagin@false\in@\tag{#1}\ifin@\tagin@true\fi}
\def\multline@#1$${\inany@true\vspace@\allowdisplaybreak@\displaybreak@
 \tagin@{#1}\iftagsleft@\DN@{\multline@l#1$$}\else
 \DN@{\multline@r#1$$}\fi\next@}
\newdimen\mwidth@
\def\rmmeasure@#1\endmultline{%
 \def\shoveleft##1{##1}\def\shoveright##1{##1}
 \setbox@ne\vbox{\Let@\halign{\setboxz@h
  {$\m@th\@lign\displaystyle{}##$}\global\mwidth@\wdz@
  \crcr#1\crcr}}}
\newdimen\mlineht@
\newif\ifzerocr@
\newif\ifonecr@
\def\lmmeasure@#1\endmultline{\global\zerocr@true\global\onecr@false
 \everycr{\noalign{\ifonecr@\global\onecr@false\fi
  \ifzerocr@\global\zerocr@false\global\onecr@true\fi}}
  \def\shoveleft##1{##1}\def\shoveright##1{##1}%
 \setbox@ne\vbox{\Let@\halign{\setboxz@h
  {$\m@th\@lign\displaystyle{}##$}\ifonecr@\global\mwidth@\wdz@
  \global\mlineht@\ht\z@\fi\crcr#1\crcr}}}
\newbox\mtagbox@
\newdimen\ltwidth@
\newdimen\rtwidth@
\def\multline@l#1$${\iftagin@\DN@{\lmultline@@#1$$}\else
 \DN@{\setbox\mtagbox@\null\ltwidth@\z@\rtwidth@\z@
  \lmultline@@@#1$$}\fi\next@}
\def\lmultline@@#1\endmultline\tag#2$${%
 \setbox\mtagbox@\hbox{\maketag@#2\maketag@}
 \lmmeasure@#1\endmultline\dimen@\mwidth@\advance\dimen@\wd\mtagbox@
 \advance\dimen@\multlinetaggap@                                            
 \ifdim\dimen@>\displaywidth\ltwidth@\z@\else\ltwidth@\wd\mtagbox@\fi       
 \lmultline@@@#1\endmultline$$}
\def\lmultline@@@{\displ@y
 \def\shoveright##1{##1\hfilneg\hskip\multlinegap@}%
 \def\shoveleft##1{\setboxz@h{$\m@th\displaystyle{}##1$}%
  \setbox@ne\hbox{$\m@th\displaystyle##1$}%
  \hfilneg
  \iftagin@
   \ifdim\ltwidth@>\z@\hskip\ltwidth@\hskip\multlinetaggap@\fi
  \else\hskip\multlinegap@\fi\hskip.5\wd@ne\hskip-.5\wdz@##1}
  \halign\bgroup\Let@\hbox to\displaywidth
   {\strut@$\m@th\displaystyle\hfil{}##\hfil$}\crcr
   \hfilneg                                                                 
   \iftagin@                                                                
    \ifdim\ltwidth@>\z@                                                     
     \box\mtagbox@\hskip\multlinetaggap@                                    
    \else
     \rlap{\vbox{\normalbaselines\hbox{\strut@\box\mtagbox@}%
     \vbox to\mlineht@{}}}\fi                                               
   \else\hskip\multlinegap@\fi}                                             
\def\multline@r#1$${\iftagin@\DN@{\rmultline@@#1$$}\else
 \DN@{\setbox\mtagbox@\null\ltwidth@\z@\rtwidth@\z@
  \rmultline@@@#1$$}\fi\next@}
\def\rmultline@@#1\endmultline\tag#2$${\ltwidth@\z@
 \setbox\mtagbox@\hbox{\maketag@#2\maketag@}%
 \rmmeasure@#1\endmultline\dimen@\mwidth@\advance\dimen@\wd\mtagbox@
 \advance\dimen@\multlinetaggap@
 \ifdim\dimen@>\displaywidth\rtwidth@\z@\else\rtwidth@\wd\mtagbox@\fi
 \rmultline@@@#1\endmultline$$}
\def\rmultline@@@{\displ@y
 \def\shoveright##1{##1\hfilneg\iftagin@\ifdim\rtwidth@>\z@
  \hskip\rtwidth@\hskip\multlinetaggap@\fi\else\hskip\multlinegap@\fi}%
 \def\shoveleft##1{\setboxz@h{$\m@th\displaystyle{}##1$}%
  \setbox@ne\hbox{$\m@th\displaystyle##1$}%
  \hfilneg\hskip\multlinegap@\hskip.5\wd@ne\hskip-.5\wdz@##1}%
 \halign\bgroup\Let@\hbox to\displaywidth
  {\strut@$\m@th\displaystyle\hfil{}##\hfil$}\crcr
 \hfilneg\hskip\multlinegap@}
\def\endmultline{\iftagsleft@\expandafter\lendmultline@\else
 \expandafter\rendmultline@\fi}
\def\lendmultline@{\hfilneg\hskip\multlinegap@\crcr\egroup}
\def\rendmultline@{\iftagin@                                                
 \ifdim\rtwidth@>\z@                                                        
  \hskip\multlinetaggap@\box\mtagbox@                                       
 \else\llap{\vtop{\normalbaselines\null\hbox{\strut@\box\mtagbox@}}}\fi     
 \else\hskip\multlinegap@\fi                                                
 \hfilneg\crcr\egroup}
\def\bmod{\mskip-\medmuskip\mkern5mu\mathbin{\fam\z@ mod}\penalty900
 \mkern5mu\mskip-\medmuskip}
\def\pmod#1{\allowbreak\ifinner\mkern8mu\else\mkern18mu\fi
 ({\fam\z@ mod}\,\,#1)}
\def\pod#1{\allowbreak\ifinner\mkern8mu\else\mkern18mu\fi(#1)}
\def\mod#1{\allowbreak\ifinner\mkern12mu\else\mkern18mu\fi{\fam\z@ mod}\,\,#1}
\message{continued fractions,}
\newcount\cfraccount@
\def\cfrac{\bgroup\bgroup\advance\cfraccount@\@ne\strut
 \iffalse{\fi\def\\{\over\displaystyle}\iffalse}\fi}
\def\lcfrac{\bgroup\bgroup\advance\cfraccount@\@ne\strut
 \iffalse{\fi\def\\{\hfill\over\displaystyle}\iffalse}\fi}
\def\rcfrac{\bgroup\bgroup\advance\cfraccount@\@ne\strut\hfill
 \iffalse{\fi\def\\{\over\displaystyle}\iffalse}\fi}
\def\gloop@#1\repeat{\gdef\body{#1}\iterate}
\def\endcfrac{\gloop@\ifnum\cfraccount@>\z@\global\advance\cfraccount@\m@ne
 \egroup\hskip-\nulldelimiterspace\egroup\repeat}
\message{compound symbols,}
\def\binrel@#1{\setboxz@h{\thinmuskip0mu
  \medmuskip\m@ne mu\thickmuskip\@ne mu$#1\m@th$}%
 \setbox@ne\hbox{\thinmuskip0mu\medmuskip\m@ne mu\thickmuskip
  \@ne mu${}#1{}\m@th$}%
 \setbox\tw@\hbox{\hskip\wd@ne\hskip-\wdz@}}
\def\overset#1\to#2{\binrel@{#2}\ifdim\wd\tw@<\z@
 \mathbin{\mathop{\kern\z@#2}\limits^{#1}}\else\ifdim\wd\tw@>\z@
 \mathrel{\mathop{\kern\z@#2}\limits^{#1}}\else
 {\mathop{\kern\z@#2}\limits^{#1}}{}\fi\fi}
\def\underset#1\to#2{\binrel@{#2}\ifdim\wd\tw@<\z@
 \mathbin{\mathop{\kern\z@#2}\limits_{#1}}\else\ifdim\wd\tw@>\z@
 \mathrel{\mathop{\kern\z@#2}\limits_{#1}}\else
 {\mathop{\kern\z@#2}\limits_{#1}}{}\fi\fi}
\def\oversetbrace#1\to#2{\overbrace{#2}^{#1}}
\def\undersetbrace#1\to#2{\underbrace{#2}_{#1}}
\def\sideset#1\and#2\to#3{%
 \setbox@ne\hbox{$\dsize{\vphantom{#3}}#1{#3}\m@th$}%
 \setbox\tw@\hbox{$\dsize{#3}#2\m@th$}%
 \hskip\wd@ne\hskip-\wd\tw@\mathop{\hskip\wd\tw@\hskip-\wd@ne
  {\vphantom{#3}}#1{#3}#2}}
\def\rightarrowfill@#1{\setboxz@h{$#1-\m@th$}\ht\z@\z@
  $#1\m@th\copy\z@\mkern-6mu\cleaders
  \hbox{$#1\mkern-2mu\box\z@\mkern-2mu$}\hfill
  \mkern-6mu\mathord\rightarrow$}
\def\leftarrowfill@#1{\setboxz@h{$#1-\m@th$}\ht\z@\z@
  $#1\m@th\mathord\leftarrow\mkern-6mu\cleaders
  \hbox{$#1\mkern-2mu\copy\z@\mkern-2mu$}\hfill
  \mkern-6mu\box\z@$}
\def\leftrightarrowfill@#1{\setboxz@h{$#1-\m@th$}\ht\z@\z@
  $#1\m@th\mathord\leftarrow\mkern-6mu\cleaders
  \hbox{$#1\mkern-2mu\box\z@\mkern-2mu$}\hfill
  \mkern-6mu\mathord\rightarrow$}
\def\overrightarrow{\mathpalette\overrightarrow@}
\def\overrightarrow@#1#2{\vbox{\ialign{##\crcr\rightarrowfill@#1\crcr
 \noalign{\kern-\ex@\nointerlineskip}$\m@th\hfil#1#2\hfil$\crcr}}}

\def\overleftarrow{\mathpalette\overleftarrow@}
\def\overleftarrow@#1#2{\vbox{\ialign{##\crcr\leftarrowfill@#1\crcr
 \noalign{\kern-\ex@\nointerlineskip}$\m@th\hfil#1#2\hfil$\crcr}}}
\def\overleftrightarrow{\mathpalette\overleftrightarrow@}
\def\overleftrightarrow@#1#2{\vbox{\ialign{##\crcr\leftrightarrowfill@#1\crcr
 \noalign{\kern-\ex@\nointerlineskip}$\m@th\hfil#1#2\hfil$\crcr}}}
\def\underrightarrow{\mathpalette\underrightarrow@}
\def\underrightarrow@#1#2{\vtop{\ialign{##\crcr$\m@th\hfil#1#2\hfil$\crcr
 \noalign{\nointerlineskip}\rightarrowfill@#1\crcr}}}

\def\underleftarrow{\mathpalette\underleftarrow@}
\def\underleftarrow@#1#2{\vtop{\ialign{##\crcr$\m@th\hfil#1#2\hfil$\crcr
 \noalign{\nointerlineskip}\leftarrowfill@#1\crcr}}}
\def\underleftrightarrow{\mathpalette\underleftrightarrow@}
\def\underleftrightarrow@#1#2{\vtop{\ialign{##\crcr$\m@th\hfil#1#2\hfil$\crcr
 \noalign{\nointerlineskip}\leftrightarrowfill@#1\crcr}}}
\message{various kinds of dots,}
\let\DOTSI\relax
\let\DOTSB\relax

\newif\ifmath@
{\uccode`7=`\\ \uccode`8=`m \uccode`9=`a \uccode`0=`t \uccode`!=`h
 \uppercase{\gdef\math@#1#2#3#4#5#6\math@{\global\math@false\ifx 7#1\ifx 8#2%
 \ifx 9#3\ifx 0#4\ifx !#5\xdef\meaning@{#6}\global\math@true\fi\fi\fi\fi\fi}}}
\newif\ifmathch@
{\uccode`7=`c \uccode`8=`h \uccode`9=`\"
 \uppercase{\gdef\mathch@#1#2#3#4#5#6\mathch@{\global\mathch@false
  \ifx 7#1\ifx 8#2\ifx 9#5\global\mathch@true\xdef\meaning@{9#6}\fi\fi\fi}}}
\newcount\classnum@
\def\getmathch@#1.#2\getmathch@{\classnum@#1 \divide\classnum@4096
 \ifcase\number\classnum@\or\or\gdef\thedots@{\dotsb@}\or
 \gdef\thedots@{\dotsb@}\fi}
\newif\ifmathbin@
{\uccode`4=`b \uccode`5=`i \uccode`6=`n
 \uppercase{\gdef\mathbin@#1#2#3{\relaxnext@
  \DNii@##1\mathbin@{\ifx\space@\next\global\mathbin@true\fi}%
 \global\mathbin@false\DN@##1\mathbin@{}%
 \ifx 4#1\ifx 5#2\ifx 6#3\DN@{\FN@\nextii@}\fi\fi\fi\next@}}}
\newif\ifmathrel@
{\uccode`4=`r \uccode`5=`e \uccode`6=`l
 \uppercase{\gdef\mathrel@#1#2#3{\relaxnext@
  \DNii@##1\mathrel@{\ifx\space@\next\global\mathrel@true\fi}%
 \global\mathrel@false\DN@##1\mathrel@{}%
 \ifx 4#1\ifx 5#2\ifx 6#3\DN@{\FN@\nextii@}\fi\fi\fi\next@}}}
\newif\ifmacro@
{\uccode`5=`m \uccode`6=`a \uccode`7=`c
 \uppercase{\gdef\macro@#1#2#3#4\macro@{\global\macro@false
  \ifx 5#1\ifx 6#2\ifx 7#3\global\macro@true
  \xdef\meaning@{\macro@@#4\macro@@}\fi\fi\fi}}}
\def\macro@@#1->#2\macro@@{#2}
\newif\ifDOTS@
\newcount\DOTSCASE@
{\uccode`6=`\\ \uccode`7=`D \uccode`8=`O \uccode`9=`T \uccode`0=`S
 \uppercase{\gdef\DOTS@#1#2#3#4#5{\global\DOTS@false\DN@##1\DOTS@{}%
  \ifx 6#1\ifx 7#2\ifx 8#3\ifx 9#4\ifx 0#5\let\next@\DOTS@@\fi\fi\fi\fi\fi
  \next@}}}
{\uccode`3=`B \uccode`4=`I \uccode`5=`X
 \uppercase{\gdef\DOTS@@#1{\relaxnext@
  \DNii@##1\DOTS@{\ifx\space@\next\global\DOTS@true\fi}%
  \DN@{\FN@\nextii@}%
  \ifx 3#1\global\DOTSCASE@\z@\else
  \ifx 4#1\global\DOTSCASE@\@ne\else
  \ifx 5#1\global\DOTSCASE@\tw@\else\DN@##1\DOTS@{}%
  \fi\fi\fi\next@}}}
\newif\ifnot@
{\uccode`5=`\\ \uccode`6=`n \uccode`7=`o \uccode`8=`t
 \uppercase{\gdef\not@#1#2#3#4{\relaxnext@
  \DNii@##1\not@{\ifx\space@\next\global\not@true\fi}%
 \global\not@false\DN@##1\not@{}%
 \ifx 5#1\ifx 6#2\ifx 7#3\ifx 8#4\DN@{\FN@\nextii@}\fi\fi\fi
 \fi\next@}}}
\newif\ifkeybin@
\def\keybin@{\keybin@true
 \ifx\next+\else\ifx\next=\else\ifx\next<\else\ifx\next>\else\ifx\next-\else
 \ifx\next*\else\ifx\next:\else\keybin@false\fi\fi\fi\fi\fi\fi\fi}
\def\dots{\RIfM@\expandafter\mdots@\else\expandafter\tdots@\fi}
\def\tdots@{\unskip\relaxnext@
 \DN@{$\m@th\mathinner{\ldotp\ldotp\ldotp}\,
   \ifx\next,\,$\else\ifx\next.\,$\else\ifx\next;\,$\else\ifx\next:\,$\else
   \ifx\next?\,$\else\ifx\next!\,$\else$ \fi\fi\fi\fi\fi\fi}%
 \ \FN@\next@}
\def\mdots@{\FN@\mdots@@}
\def\mdots@@{\gdef\thedots@{\dotso@}
 \ifx\next\boldkey\gdef\thedots@\boldkey{\boldkeydots@}\else                
 \ifx\next\boldsymbol\gdef\thedots@\boldsymbol{\boldsymboldots@}\else       
 \ifx,\next\gdef\thedots@{\dotsc}
 \else\ifx\not\next\gdef\thedots@{\dotsb@}
 \else\keybin@
 \ifkeybin@\gdef\thedots@{\dotsb@}
 \else\xdef\meaning@{\meaning\next..........}\xdef\meaning@@{\meaning@}
  \expandafter\math@\meaning@\math@
  \ifmath@
   \expandafter\mathch@\meaning@\mathch@
   \ifmathch@\expandafter\getmathch@\meaning@\getmathch@\fi                 
  \else\expandafter\macro@\meaning@@\macro@                                 
  \ifmacro@                                                                
   \expandafter\not@\meaning@\not@\ifnot@\gdef\thedots@{\dotsb@}
  \else\expandafter\DOTS@\meaning@\DOTS@
  \ifDOTS@
   \ifcase\number\DOTSCASE@\gdef\thedots@{\dotsb@}%
    \or\gdef\thedots@{\dotsi}\else\fi                                      
  \else\expandafter\math@\meaning@\math@                                   
  \ifmath@\expandafter\mathbin@\meaning@\mathbin@
  \ifmathbin@\gdef\thedots@{\dotsb@}
  \else\expandafter\mathrel@\meaning@\mathrel@
  \ifmathrel@\gdef\thedots@{\dotsb@}
  \fi\fi\fi\fi\fi\fi\fi\fi\fi\fi\fi\fi
 \thedots@}
\def\plainldots@{\mathinner{\ldotp\ldotp\ldotp}}
\def\plaincdots@{\mathinner{\cdotp\cdotp\cdotp}}
\def\dotsi{\!\plaincdots@}
\let\dotsb@\plaincdots@
\newif\ifextra@
\newif\ifrightdelim@
\def\rightdelim@{\global\rightdelim@true                                    
 \ifx\next)\else                                                            
 \ifx\next]\else
 \ifx\next\rbrack\else
 \ifx\next\}\else
 \ifx\next\rbrace\else
 \ifx\next\rangle\else
 \ifx\next\rceil\else
 \ifx\next\rfloor\else
 \ifx\next\rgroup\else
 \ifx\next\rmoustache\else
 \ifx\next\right\else
 \ifx\next\bigr\else
 \ifx\next\biggr\else
 \ifx\next\Bigr\else                                                        
 \ifx\next\Biggr\else\global\rightdelim@false
 \fi\fi\fi\fi\fi\fi\fi\fi\fi\fi\fi\fi\fi\fi\fi}
\def\extra@{%
 \global\extra@false\rightdelim@\ifrightdelim@\global\extra@true            
 \else\ifx\next$\global\extra@true                                          
 \else\xdef\meaning@{\meaning\next..........}
 \expandafter\macro@\meaning@\macro@\ifmacro@                               
 \expandafter\DOTS@\meaning@\DOTS@
 \ifDOTS@
 \ifnum\DOTSCASE@=\tw@\global\extra@true                                    
 \fi\fi\fi\fi\fi}
\newif\ifbold@
\def\dotso@{\relaxnext@
 \ifbold@
  \let\next\delayed@
  \DNii@{\extra@\plainldots@\ifextra@\,\fi}%
 \else
  \DNii@{\DN@{\extra@\plainldots@\ifextra@\,\fi}\FN@\next@}%
 \fi
 \nextii@}
\def\extrap@#1{%
 \ifx\next,\DN@{#1\,}\else
 \ifx\next;\DN@{#1\,}\else
 \ifx\next.\DN@{#1\,}\else\extra@
 \ifextra@\DN@{#1\,}\else
 \let\next@#1\fi\fi\fi\fi\next@}
\def\ldots{\DN@{\extrap@\plainldots@}%
 \FN@\next@}
\def\cdots{\DN@{\extrap@\plaincdots@}%
 \FN@\next@}

\def\dotsc{\relaxnext@
 \DN@{\ifx\next;\plainldots@\,\else
  \ifx\next.\plainldots@\,\else\extra@\plainldots@
  \ifextra@\,\fi\fi\fi}%
 \FN@\next@}
\def\cdot{\mathchar"2201 }

\message{special superscripts,}
\def\dddot#1{{\mathop{#1}\limits^{\vbox to-1.4\ex@{\kern-\tw@\ex@
 \hbox{\rm...}\vss}}}}
\def\ddddot#1{{\mathop{#1}\limits^{\vbox to-1.4\ex@{\kern-\tw@\ex@
 \hbox{\rm....}\vss}}}}
\def\sphat{^{\mathchoice{}{}%
 {\,\,\botsmash{\hbox{\lower4\ex@\hbox{$\m@th\widehat{\null}$}}}}%
 {\,\botsmash{\hbox{\lower3\ex@\hbox{$\m@th\hat{\null}$}}}}}}

\def\spacute{^{\!\botsmash{\hbox{\lower\@ne ex\hbox{\'{}}}}}}
\def\spgrave{^{\mathchoice{}{}{}{\!}%
 \botsmash{\hbox{\lower\@ne ex\hbox{\`{}}}}}}
\def\spdot{^{\hbox{\raise\ex@\hbox{\rm.}}}}
\def\spddot{^{\hbox{\raise\ex@\hbox{\rm..}}}}
\def\spdddot{^{\hbox{\raise\ex@\hbox{\rm...}}}}
\def\spddddot{^{\hbox{\raise\ex@\hbox{\rm....}}}}
\def\spbreve{^{\!\botsmash{\hbox{\lower4\ex@\hbox{\u{}}}}}}

\message{\string\text,}
\def\textonlyfont@#1#2{\def#1{\RIfM@
 \Err@{Use \string#1\space only in text}\else#2\fi}}
\textonlyfont@\rm\tenrm
\textonlyfont@\it\tenit
\textonlyfont@\sl\tensl
\textonlyfont@\bf\tenbf
\def\oldnos#1{\RIfM@{\mathcode`\,="013B \fam\@ne#1}\else
 \leavevmode\hbox{$\m@th\mathcode`\,="013B \fam\@ne#1$}\fi}
\def\text{\RIfM@\expandafter\text@\else\expandafter\text@@\fi}
\def\text@@#1{\leavevmode\hbox{#1}}
\def\mathhexbox@#1#2#3{\text{$\m@th\mathchar"#1#2#3$}}
\def\dag{{\mathhexbox@279}}
\def\ddag{{\mathhexbox@27A}}
\def\S{{\mathhexbox@278}}
\def\P{{\mathhexbox@27B}}
\newif\iffirstchoice@
\firstchoice@true
\def\text@#1{\mathchoice
 {\hbox{\everymath{\displaystyle}\def\textfonti{\the\textfont\@ne}%
  \def\textfontii{\the\textfont\tw@}\textdef@@ T#1}}
 {\hbox{\firstchoice@false
  \everymath{\textstyle}\def\textfonti{\the\textfont\@ne}%
  \def\textfontii{\the\textfont\tw@}\textdef@@ T#1}}
 {\hbox{\firstchoice@false
  \everymath{\scriptstyle}\def\textfonti{\the\scriptfont\@ne}%
  \def\textfontii{\the\scriptfont\tw@}\textdef@@ S\rm#1}}
 {\hbox{\firstchoice@false
  \everymath{\scriptscriptstyle}\def\textfonti
  {\the\scriptscriptfont\@ne}%
  \def\textfontii{\the\scriptscriptfont\tw@}\textdef@@ s\rm#1}}}
\def\textdef@@#1{\textdef@#1\rm\textdef@#1\bf\textdef@#1\sl\textdef@#1\it}
\def\rmfam{0}
\def\textdef@#1#2{%
 \DN@{\csname\expandafter\eat@\string#2fam\endcsname}%
 \if S#1\edef#2{\the\scriptfont\next@\relax}%
 \else\if s#1\edef#2{\the\scriptscriptfont\next@\relax}%
 \else\edef#2{\the\textfont\next@\relax}\fi\fi}
\scriptfont\itfam\tenit \scriptscriptfont\itfam\tenit
\scriptfont\slfam\tensl \scriptscriptfont\slfam\tensl
\newif\iftopfolded@
\newif\ifbotfolded@
\def\topfoldedtext{\topfolded@true\botfolded@false\foldedtext@}
\def\botfoldedtext{\botfolded@true\topfolded@false\foldedtext@}
\def\foldedtext{\topfolded@false\botfolded@false\foldedtext@}
\Invalid@\foldedwidth
\def\foldedtext@{\relaxnext@
 \DN@{\ifx\next\foldedwidth\let\next@\nextii@\else
  \DN@{\nextii@\foldedwidth{.3\hsize}}\fi\next@}%
 \DNii@\foldedwidth##1##2{\setbox\z@\vbox
  {\normalbaselines\hsize##1\relax
  \tolerance1600 \noindent\ignorespaces##2}\ifbotfolded@\boxz@\else
  \iftopfolded@\vtop{\unvbox\z@}\else\vcenter{\boxz@}\fi\fi}%
 \FN@\next@}
\message{math font commands,}
\def\bold{\RIfM@\expandafter\bold@\else
 \expandafter\nonmatherr@\expandafter\bold\fi}
\def\bold@#1{{\bold@@{#1}}}
\def\bold@@#1{\fam\bffam\relax#1}
\def\slanted{\RIfM@\expandafter\slanted@\else
 \expandafter\nonmatherr@\expandafter\slanted\fi}
\def\slanted@#1{{\slanted@@{#1}}}
\def\slanted@@#1{\fam\slfam\relax#1}
\def\roman{\RIfM@\expandafter\roman@\else
 \expandafter\nonmatherr@\expandafter\roman\fi}
\def\roman@#1{{\roman@@{#1}}}
\def\roman@@#1{\fam\rmfam\relax#1}
\def\italic{\RIfM@\expandafter\italic@\else
 \expandafter\nonmatherr@\expandafter\italic\fi}
\def\italic@#1{{\italic@@{#1}}}
\def\italic@@#1{\fam\itfam\relax#1}
\def\Cal{\RIfM@\expandafter\Cal@\else
 \expandafter\nonmatherr@\expandafter\Cal\fi}
\def\Cal@#1{{\Cal@@{#1}}}
\def\Cal@@#1{\noaccents@\fam\tw@#1}
\mathchardef\Gamma="0000
\mathchardef\Delta="0001
\mathchardef\Theta="0002
\mathchardef\Lambda="0003
\mathchardef\Xi="0004
\mathchardef\Pi="0005
\mathchardef\Sigma="0006
\mathchardef\Upsilon="0007
\mathchardef\Phi="0008
\mathchardef\Psi="0009
\mathchardef\Omega="000A
\mathchardef\varGamma="0100
\mathchardef\varDelta="0101
\mathchardef\varTheta="0102
\mathchardef\varLambda="0103
\mathchardef\varXi="0104
\mathchardef\varPi="0105
\mathchardef\varSigma="0106
\mathchardef\varUpsilon="0107
\mathchardef\varPhi="0108
\mathchardef\varPsi="0109
\mathchardef\varOmega="010A
\let\alloc@@\alloc@
\def\hexnumber@#1{\ifcase#1 0\or 1\or 2\or 3\or 4\or 5\or 6\or 7\or 8\or
 9\or A\or B\or C\or D\or E\or F\fi}
\def\loadmsam{%
 \font@\tenmsa=msam10
 \font@\sevenmsa=msam7
 \font@\fivemsa=msam5
 \alloc@@8\fam\chardef\sixt@@n\msafam
 \textfont\msafam=\tenmsa
 \scriptfont\msafam=\sevenmsa
 \scriptscriptfont\msafam=\fivemsa
 \edef\next{\hexnumber@\msafam}%
 \mathchardef\dabar@"0\next39
 \edef\dashrightarrow{\mathrel{\dabar@\dabar@\mathchar"0\next4B}}%
 \edef\dashleftarrow{\mathrel{\mathchar"0\next4C\dabar@\dabar@}}%
 \let\dasharrow\dashrightarrow
 \edef\ulcorner{\delimiter"4\next70\next70 }%
 \edef\urcorner{\delimiter"5\next71\next71 }%
 \edef\llcorner{\delimiter"4\next78\next78 }%
 \edef\lrcorner{\delimiter"5\next79\next79 }%
 \edef\yen{{\noexpand\mathhexbox@\next55}}%
 \edef\checkmark{{\noexpand\mathhexbox@\next58}}%
 \edef\circledR{{\noexpand\mathhexbox@\next72}}%
 \edef\maltese{{\noexpand\mathhexbox@\next7A}}%
 \global\let\loadmsam\empty}%
\def\loadmsbm{%
 \font@\tenmsb=msbm10 \font@\sevenmsb=msbm7 \font@\fivemsb=msbm5
 \alloc@@8\fam\chardef\sixt@@n\msbfam
 \textfont\msbfam=\tenmsb
 \scriptfont\msbfam=\sevenmsb \scriptscriptfont\msbfam=\fivemsb
 \global\let\loadmsbm\empty
 }
\def\widehat#1{\ifx\undefined\msbfam \DN@{362}%
  \else \setboxz@h{$\m@th#1$}%
    \edef\next@{\ifdim\wdz@>\tw@ em%
        \hexnumber@\msbfam 5B%
      \else 362\fi}\fi
  \mathaccent"0\next@{#1}}
\def\widetilde#1{\ifx\undefined\msbfam \DN@{365}%
  \else \setboxz@h{$\m@th#1$}%
    \edef\next@{\ifdim\wdz@>\tw@ em%
        \hexnumber@\msbfam 5D%
      \else 365\fi}\fi
  \mathaccent"0\next@{#1}}
\message{\string\newsymbol,}
\def\newsymbol#1#2#3#4#5{\define#1{}%
  \count@#2\relax \advance\count@\m@ne 
 \ifcase\count@
   \ifx\undefined\msafam\loadmsam\fi \let\next@\msafam
 \or \ifx\undefined\msbfam\loadmsbm\fi \let\next@\msbfam
 \else  \Err@{\Invalid@@\string\newsymbol}\let\next@\tw@\fi
 \mathchardef#1="#3\hexnumber@\next@#4#5\space}

\def\Bbb{\RIfM@\expandafter\Bbb@\else
 \expandafter\nonmatherr@\expandafter\Bbb\fi}
\def\Bbb@#1{{\Bbb@@{#1}}}
\def\Bbb@@#1{\noaccents@\fam\msbfam\relax#1}
\message{bold Greek and bold symbols,}
\def\loadbold{%
 \font@\tencmmib=cmmib10 \font@\sevencmmib=cmmib7 \font@\fivecmmib=cmmib5
 \skewchar\tencmmib'177 \skewchar\sevencmmib'177 \skewchar\fivecmmib'177
 \alloc@@8\fam\chardef\sixt@@n\cmmibfam
 \textfont\cmmibfam\tencmmib
 \scriptfont\cmmibfam\sevencmmib \scriptscriptfont\cmmibfam\fivecmmib
 \font@\tencmbsy=cmbsy10 \font@\sevencmbsy=cmbsy7 \font@\fivecmbsy=cmbsy5
 \skewchar\tencmbsy'60 \skewchar\sevencmbsy'60 \skewchar\fivecmbsy'60
 \alloc@@8\fam\chardef\sixt@@n\cmbsyfam
 \textfont\cmbsyfam\tencmbsy
 \scriptfont\cmbsyfam\sevencmbsy \scriptscriptfont\cmbsyfam\fivecmbsy
 \let\loadbold\empty
}
\def\boldnotloaded#1{\Err@{\ifcase#1\or First\else Second\fi
       bold symbol font not loaded}}
\def\mathchari@#1#2#3{\ifx\undefined\cmmibfam
    \boldnotloaded@\@ne
  \else\mathchar"#1\hexnumber@\cmmibfam#2#3\space \fi}
\def\mathcharii@#1#2#3{\ifx\undefined\cmbsyfam
    \boldnotloaded\tw@
  \else \mathchar"#1\hexnumber@\cmbsyfam#2#3\space\fi}
\edef\bffam@{\hexnumber@\bffam}
\def\boldkey#1{\ifcat\noexpand#1A%
  \ifx\undefined\cmmibfam \boldnotloaded\@ne
  \else {\fam\cmmibfam#1}\fi
 \else
 \ifx#1!\mathchar"5\bffam@21 \else
 \ifx#1(\mathchar"4\bffam@28 \else\ifx#1)\mathchar"5\bffam@29 \else
 \ifx#1+\mathchar"2\bffam@2B \else\ifx#1:\mathchar"3\bffam@3A \else
 \ifx#1;\mathchar"6\bffam@3B \else\ifx#1=\mathchar"3\bffam@3D \else
 \ifx#1?\mathchar"5\bffam@3F \else\ifx#1[\mathchar"4\bffam@5B \else
 \ifx#1]\mathchar"5\bffam@5D \else
 \ifx#1,\mathchari@63B \else
 \ifx#1-\mathcharii@200 \else
 \ifx#1.\mathchari@03A \else
 \ifx#1/\mathchari@03D \else
 \ifx#1<\mathchari@33C \else
 \ifx#1>\mathchari@33E \else
 \ifx#1*\mathcharii@203 \else
 \ifx#1|\mathcharii@06A \else
 \ifx#10\bold0\else\ifx#11\bold1\else\ifx#12\bold2\else\ifx#13\bold3\else
 \ifx#14\bold4\else\ifx#15\bold5\else\ifx#16\bold6\else\ifx#17\bold7\else
 \ifx#18\bold8\else\ifx#19\bold9\else
  \Err@{\string\boldkey\space can't be used with #1}%
 \fi\fi\fi\fi\fi\fi\fi\fi\fi\fi\fi\fi\fi\fi\fi
 \fi\fi\fi\fi\fi\fi\fi\fi\fi\fi\fi\fi\fi\fi}
\def\boldsymbol#1{%
 \DN@{\Err@{You can't use \string\boldsymbol\space with \string#1}#1}%
 \ifcat\noexpand#1A%
   \let\next@\relax
   \ifx\undefined\cmmibfam \boldnotloaded\@ne
   \else {\fam\cmmibfam#1}\fi
 \else
  \xdef\meaning@{\meaning#1.........}%
  \expandafter\math@\meaning@\math@
  \ifmath@
   \expandafter\mathch@\meaning@\mathch@
   \ifmathch@
    \expandafter\boldsymbol@@\meaning@\boldsymbol@@
   \fi
  \else
   \expandafter\macro@\meaning@\macro@
   \expandafter\delim@\meaning@\delim@
   \ifdelim@
    \expandafter\delim@@\meaning@\delim@@
   \else
    \boldsymbol@{#1}%
   \fi
  \fi
 \fi
 \next@}
\def\mathhexboxii@#1#2{\ifx\undefined\cmbsyfam
    \boldnotloaded\tw@
  \else \mathhexbox@{\hexnumber@\cmbsyfam}{#1}{#2}\fi}
\def\boldsymbol@#1{\let\next@\relax\let\next#1%
 \ifx\next\cdot\mathcharii@201 \else
 \ifx\next\prime{{\null\mathcharii@030 \null}}\else
 \ifx\next\lbrack\mathchar"4\bffam@5B \else
 \ifx\next\rbrack\mathchar"5\bffam@5D \else
 \ifx\next\{\mathcharii@466 \else
 \ifx\next\lbrace\mathcharii@466 \else
 \ifx\next\}\mathcharii@567 \else
 \ifx\next\rbrace\mathcharii@567 \else
 \ifx\next\surd{{\mathcharii@170}}\else
 \ifx\next\S{{\mathhexboxii@78}}\else
 \ifx\next\P{{\mathhexboxii@7B}}\else
 \ifx\next\dag{{\mathhexboxii@79}}\else
 \ifx\next\ddag{{\mathhexboxii@7A}}\else
 \DN@{\Err@{You can't use \string\boldsymbol\space with \string#1}#1}%
 \fi\fi\fi\fi\fi\fi\fi\fi\fi\fi\fi\fi\fi}
\def\boldsymbol@@#1.#2\boldsymbol@@{\classnum@#1 \count@@@\classnum@        
 \divide\classnum@4096 \count@\classnum@                                    
 \multiply\count@4096 \advance\count@@@-\count@ \count@@\count@@@           
 \divide\count@@@\@cclvi \count@\count@@                                    
 \multiply\count@@@\@cclvi \advance\count@@-\count@@@                       
 \divide\count@@@\@cclvi                                                    
 \multiply\classnum@4096 \advance\classnum@\count@@                         
 \ifnum\count@@@=\z@                                                        
  \count@"\bffam@ \multiply\count@\@cclvi
  \advance\classnum@\count@
  \DN@{\mathchar\number\classnum@}%
 \else
  \ifnum\count@@@=\@ne                                                      
   \ifx\undefined\cmmibfam \DN@{\boldnotloaded\@ne}%
   \else \count@\cmmibfam \multiply\count@\@cclvi
     \advance\classnum@\count@
     \DN@{\mathchar\number\classnum@}\fi
  \else
   \ifnum\count@@@=\tw@                                                    
     \ifx\undefined\cmbsyfam
       \DN@{\boldnotloaded\tw@}%
     \else
       \count@\cmbsyfam \multiply\count@\@cclvi
       \advance\classnum@\count@
       \DN@{\mathchar\number\classnum@}%
     \fi
  \fi
 \fi
\fi}
\newif\ifdelim@
\newcount\delimcount@
{\uccode`6=`\\ \uccode`7=`d \uccode`8=`e \uccode`9=`l
 \uppercase{\gdef\delim@#1#2#3#4#5\delim@
  {\delim@false\ifx 6#1\ifx 7#2\ifx 8#3\ifx 9#4\delim@true
   \xdef\meaning@{#5}\fi\fi\fi\fi}}}
\def\delim@@#1"#2#3#4#5#6\delim@@{\if#32%
\let\next@\relax
 \ifx\undefined\cmbsyfam \boldnotloaded\@ne
 \else \mathcharii@#2#4#5\space \fi\fi}
\def\vert{\delimiter"026A30C }
\def\Vert{\delimiter"026B30D }
\let\|\Vert
\def\backslash{\delimiter"026E30F }
\def\boldkeydots@#1{\bold@true\let\next=#1\let\delayed@=#1\mdots@@
 \boldkey#1\bold@false}  
\def\boldsymboldots@#1{\bold@true\let\next#1\let\delayed@#1\mdots@@
 \boldsymbol#1\bold@false}
\message{Euler fonts,}

\def\frak{\mathfont@\frak}

\def\loadmathfont#1{%
   \expandafter\font@\csname ten#1\endcsname=#110
   \expandafter\font@\csname seven#1\endcsname=#17
   \expandafter\font@\csname five#1\endcsname=#15
   \edef\next{\noexpand\alloc@@8\fam\chardef\sixt@@n
     \expandafter\noexpand\csname#1fam\endcsname}%
   \next
   \textfont\csname#1fam\endcsname \csname ten#1\endcsname
   \scriptfont\csname#1fam\endcsname \csname seven#1\endcsname
   \scriptscriptfont\csname#1fam\endcsname \csname five#1\endcsname
   \expandafter\def\csname #1\expandafter\endcsname\expandafter{%
      \expandafter\mathfont@\csname#1\endcsname}%
 \expandafter\gdef\csname load#1\endcsname{}%
}
\def\mathfont@#1{\RIfM@\expandafter\mathfont@@\expandafter#1\else
  \expandafter\nonmatherr@\expandafter#1\fi}
\def\mathfont@@#1#2{{\mathfont@@@#1{#2}}}
\def\mathfont@@@#1#2{\noaccents@
   \fam\csname\expandafter\eat@\string#1fam\endcsname
   \relax#2}
\message{math accents,}
\def\accentclass@{7}
\def\noaccents@{\def\accentclass@{0}}
\def\makeacc@#1#2{\def#1{\mathaccent"\accentclass@#2 }}
\makeacc@\hat{05E}
\makeacc@\check{014}
\makeacc@\tilde{07E}
\makeacc@\acute{013}
\makeacc@\grave{012}
\makeacc@\dot{05F}
\makeacc@\ddot{07F}
\makeacc@\breve{015}
\makeacc@\bar{016}

\newcount\skewcharcount@
\newcount\familycount@
\def\theskewchar@{\familycount@\@ne
 \global\skewcharcount@\the\skewchar\textfont\@ne                           
 \ifnum\fam>\m@ne\ifnum\fam<16
  \global\familycount@\the\fam\relax
  \global\skewcharcount@\the\skewchar\textfont\the\fam\relax\fi\fi          
 \ifnum\skewcharcount@>\m@ne
  \ifnum\skewcharcount@<128
  \multiply\familycount@256
  \global\advance\skewcharcount@\familycount@
  \global\advance\skewcharcount@28672
  \mathchar\skewcharcount@\else
  \global\skewcharcount@\m@ne\fi\else
 \global\skewcharcount@\m@ne\fi}                                            
\newcount\pointcount@
\def\getpoints@#1.#2\getpoints@{\pointcount@#1 }
\newdimen\accentdimen@
\newcount\accentmu@
\def\dimentomu@{\multiply\accentdimen@ 100
 \expandafter\getpoints@\the\accentdimen@\getpoints@
 \multiply\pointcount@18
 \divide\pointcount@\@m
 \global\accentmu@\pointcount@}
\def\Makeacc@#1#2{\def#1{\RIfM@\DN@{\mathaccent@
 {"\accentclass@#2 }}\else\DN@{\nonmatherr@{#1}}\fi\next@}}
\def\unbracefonts@{\let\Cal@\Cal@@\let\roman@\roman@@\let\bold@\bold@@
 \let\slanted@\slanted@@}
\def\mathaccent@#1#2{\ifnum\fam=\m@ne\xdef\thefam@{1}\else
 \xdef\thefam@{\the\fam}\fi                                                 
 \accentdimen@\z@                                                           
 \setboxz@h{\unbracefonts@$\m@th\fam\thefam@\relax#2$}
 \ifdim\accentdimen@=\z@\DN@{\mathaccent#1{#2}}
  \setbox@ne\hbox{\unbracefonts@$\m@th\fam\thefam@\relax#2\theskewchar@$}
  \setbox\tw@\hbox{$\m@th\ifnum\skewcharcount@=\m@ne\else
   \mathchar\skewcharcount@\fi$}
  \global\accentdimen@\wd@ne\global\advance\accentdimen@-\wdz@
  \global\advance\accentdimen@-\wd\tw@                                     
  \global\multiply\accentdimen@\tw@
  \dimentomu@\global\advance\accentmu@\@ne                                 
 \else\DN@{{\mathaccent#1{#2\mkern\accentmu@ mu}%
    \mkern-\accentmu@ mu}{}}\fi                                             
 \next@}\Makeacc@\Hat{05E}
\Makeacc@\Check{014}
\Makeacc@\Tilde{07E}
\Makeacc@\Acute{013}
\Makeacc@\Grave{012}
\Makeacc@\Dot{05F}
\Makeacc@\Ddot{07F}
\Makeacc@\Breve{015}
\Makeacc@\Bar{016}
\def\Vec{\RIfM@\DN@{\mathaccent@{"017E }}\else
 \DN@{\nonmatherr@\Vec}\fi\next@}
\def\accentedsymbol#1#2{\csname newbox\expandafter\endcsname
  \csname\expandafter\eat@\string#1@box\endcsname
 \expandafter\setbox\csname\expandafter\eat@
  \string#1@box\endcsname\hbox{$\m@th#2$}\define
  #1{\copy\csname\expandafter\eat@\string#1@box\endcsname{}}}
\message{roots,}
\def\sqrt#1{\radical"270370 {#1}}
\let\underline@\underline
\let\overline@\overline
\def\underline#1{\underline@{#1}}
\def\overline#1{\overline@{#1}}
\Invalid@\leftroot
\Invalid@\uproot
\newcount\uproot@
\newcount\leftroot@
\def\root{\relaxnext@
  \DN@{\ifx\next\uproot\let\next@\nextii@\else
   \ifx\next\leftroot\let\next@\nextiii@\else
   \let\next@\plainroot@\fi\fi\next@}%
  \DNii@\uproot##1{\uproot@##1\relax\FN@\nextiv@}%
  \def\nextiv@{\ifx\next\space@\DN@. {\FN@\nextv@}\else
   \DN@.{\FN@\nextv@}\fi\next@.}%
  \def\nextv@{\ifx\next\leftroot\let\next@\nextvi@\else
   \let\next@\plainroot@\fi\next@}%
  \def\nextvi@\leftroot##1{\leftroot@##1\relax\plainroot@}%
   \def\nextiii@\leftroot##1{\leftroot@##1\relax\FN@\nextvii@}%
  \def\nextvii@{\ifx\next\space@
   \DN@. {\FN@\nextviii@}\else
   \DN@.{\FN@\nextviii@}\fi\next@.}%
  \def\nextviii@{\ifx\next\uproot\let\next@\nextix@\else
   \let\next@\plainroot@\fi\next@}%
  \def\nextix@\uproot##1{\uproot@##1\relax\plainroot@}%
  \bgroup\uproot@\z@\leftroot@\z@\FN@\next@}
\def\plainroot@#1\of#2{\setbox\rootbox\hbox{$\m@th\scriptscriptstyle{#1}$}%
 \mathchoice{\r@@t\displaystyle{#2}}{\r@@t\textstyle{#2}}
 {\r@@t\scriptstyle{#2}}{\r@@t\scriptscriptstyle{#2}}\egroup}
\def\r@@t#1#2{\setboxz@h{$\m@th#1\sqrt{#2}$}%
 \dimen@\ht\z@\advance\dimen@-\dp\z@
 \setbox@ne\hbox{$\m@th#1\mskip\uproot@ mu$}\advance\dimen@ 1.667\wd@ne
 \mkern-\leftroot@ mu\mkern5mu\raise.6\dimen@\copy\rootbox
 \mkern-10mu\mkern\leftroot@ mu\boxz@}
\def\boxed#1{\setboxz@h{$\m@th\displaystyle{#1}$}\dimen@.4\ex@
 \advance\dimen@3\ex@\advance\dimen@\dp\z@
 \hbox{\lower\dimen@\hbox{%
 \vbox{\hrule height.4\ex@
 \hbox{\vrule width.4\ex@\hskip3\ex@\vbox{\vskip3\ex@\boxz@\vskip3\ex@}%
 \hskip3\ex@\vrule width.4\ex@}\hrule height.4\ex@}%
 }}}
\message{commutative diagrams,}
\let\ampersand@\relax
\newdimen\minaw@
\minaw@11.11128\ex@
\newdimen\minCDaw@
\minCDaw@2.5pc
\def\minCDarrowwidth#1{\RIfMIfI@\onlydmatherr@\minCDarrowwidth
 \else\minCDaw@#1\relax\fi\else\onlydmatherr@\minCDarrowwidth\fi}
\newif\ifCD@
\def\CD{\bgroup\vspace@\relax\let\ampersand@&\iffalse}\fi
 \CD@true\vcenter\bgroup\Let@\tabskip\z@skip\baselineskip20\ex@
 \lineskip3\ex@\lineskiplimit3\ex@\halign\bgroup
 &\hfill$\m@th##$\hfill\crcr}
\def\endCD{\crcr\egroup\egroup\egroup}
\newdimen\bigaw@
\atdef@>#1>#2>{\ampersand@                                                  
 \setboxz@h{$\m@th\ssize\;{#1}\;\;$}
 \setbox@ne\hbox{$\m@th\ssize\;{#2}\;\;$}
 \setbox\tw@\hbox{$\m@th#2$}
 \ifCD@\global\bigaw@\minCDaw@\else\global\bigaw@\minaw@\fi                 
 \ifdim\wdz@>\bigaw@\global\bigaw@\wdz@\fi
 \ifdim\wd@ne>\bigaw@\global\bigaw@\wd@ne\fi                                
 \ifCD@\enskip\fi                                                           
 \ifdim\wd\tw@>\z@
  \mathrel{\mathop{\hbox to\bigaw@{\rightarrowfill@\displaystyle}}%
    \limits^{#1}_{#2}}
 \else\mathrel{\mathop{\hbox to\bigaw@{\rightarrowfill@\displaystyle}}%
    \limits^{#1}}\fi                                                        
 \ifCD@\enskip\fi                                                          
 \ampersand@}                                                              
\atdef@<#1<#2<{\ampersand@\setboxz@h{$\m@th\ssize\;\;{#1}\;$}%
 \setbox@ne\hbox{$\m@th\ssize\;\;{#2}\;$}\setbox\tw@\hbox{$\m@th#2$}%
 \ifCD@\global\bigaw@\minCDaw@\else\global\bigaw@\minaw@\fi
 \ifdim\wdz@>\bigaw@\global\bigaw@\wdz@\fi
 \ifdim\wd@ne>\bigaw@\global\bigaw@\wd@ne\fi
 \ifCD@\enskip\fi
 \ifdim\wd\tw@>\z@
  \mathrel{\mathop{\hbox to\bigaw@{\leftarrowfill@\displaystyle}}%
       \limits^{#1}_{#2}}\else
  \mathrel{\mathop{\hbox to\bigaw@{\leftarrowfill@\displaystyle}}%
       \limits^{#1}}\fi
 \ifCD@\enskip\fi\ampersand@}
\begingroup
 \catcode`\~=\active \lccode`\~=`\@
 \lowercase{%
  \global\atdef@)#1)#2){~>#1>#2>}
  \global\atdef@(#1(#2({~<#1<#2<}}
\endgroup
\atdef@ A#1A#2A{\llap{$\m@th\vcenter{\hbox
 {$\ssize#1$}}$}\Big\uparrow\rlap{$\m@th\vcenter{\hbox{$\ssize#2$}}$}&&}
\atdef@ V#1V#2V{\llap{$\m@th\vcenter{\hbox
 {$\ssize#1$}}$}\Big\downarrow\rlap{$\m@th\vcenter{\hbox{$\ssize#2$}}$}&&}
\atdef@={&\enskip\mathrel
 {\vbox{\hrule width\minCDaw@\vskip3\ex@\hrule width
 \minCDaw@}}\enskip&}
\atdef@|{\Big\Vert&&}
\atdef@\vert{\Big\Vert&&}
\def\pretend#1\haswidth#2{\setboxz@h{$\m@th\scriptstyle{#2}$}\hbox
 to\wdz@{\hfill$\m@th\scriptstyle{#1}$\hfill}}
\message{poor man's bold,}
\def\pmb{\RIfM@\expandafter\mathpalette\expandafter\pmb@\else
 \expandafter\pmb@@\fi}
\def\pmb@@#1{\leavevmode\setboxz@h{#1}%
   \dimen@-\wdz@
   \kern-.5\ex@\copy\z@
   \kern\dimen@\kern.25\ex@\raise.4\ex@\copy\z@
   \kern\dimen@\kern.25\ex@\box\z@
}
\def\binrel@@#1{\ifdim\wd2<\z@\mathbin{#1}\else\ifdim\wd\tw@>\z@
 \mathrel{#1}\else{#1}\fi\fi}
\newdimen\pmbraise@
\def\pmb@#1#2{\setbox\thr@@\hbox{$\m@th#1{#2}$}%
 \setbox4\hbox{$\m@th#1\mkern.5mu$}\pmbraise@\wd4\relax
 \binrel@{#2}%
 \dimen@-\wd\thr@@
   \binrel@@{%
   \mkern-.8mu\copy\thr@@
   \kern\dimen@\mkern.4mu\raise\pmbraise@\copy\thr@@
   \kern\dimen@\mkern.4mu\box\thr@@
}}
\def\documentstyle#1{\W@{}\input #1.sty\relax}
\message{syntax check,}
\font\dummyft@=dummy
\fontdimen1 \dummyft@=\z@
\fontdimen2 \dummyft@=\z@
\fontdimen3 \dummyft@=\z@
\fontdimen4 \dummyft@=\z@
\fontdimen5 \dummyft@=\z@
\fontdimen6 \dummyft@=\z@
\fontdimen7 \dummyft@=\z@
\fontdimen8 \dummyft@=\z@
\fontdimen9 \dummyft@=\z@
\fontdimen10 \dummyft@=\z@
\fontdimen11 \dummyft@=\z@
\fontdimen12 \dummyft@=\z@
\fontdimen13 \dummyft@=\z@
\fontdimen14 \dummyft@=\z@
\fontdimen15 \dummyft@=\z@
\fontdimen16 \dummyft@=\z@
\fontdimen17 \dummyft@=\z@
\fontdimen18 \dummyft@=\z@
\fontdimen19 \dummyft@=\z@
\fontdimen20 \dummyft@=\z@
\fontdimen21 \dummyft@=\z@
\fontdimen22 \dummyft@=\z@
\def\fontlist@{\\{\tenrm}\\{\sevenrm}\\{\fiverm}\\{\teni}\\{\seveni}%
 \\{\fivei}\\{\tensy}\\{\sevensy}\\{\fivesy}\\{\tenex}\\{\tenbf}\\{\sevenbf}%
 \\{\fivebf}\\{\tensl}\\{\tenit}}
\def\font@#1=#2 {\rightappend@#1\to\fontlist@\font#1=#2 }
\def\dodummy@{{\def\\##1{\global\let##1\dummyft@}\fontlist@}}
\def\nopages@{\output{\setbox\z@\box\@cclv \deadcycles\z@}%
 \alloc@5\toks\toksdef\@cclvi\output}
\let\galleys\nopages@
\newif\ifsyntax@
\newcount\countxviii@
\def\syntax{\syntax@true\dodummy@\countxviii@\count18
 \loop\ifnum\countxviii@>\m@ne\textfont\countxviii@=\dummyft@
 \scriptfont\countxviii@=\dummyft@\scriptscriptfont\countxviii@=\dummyft@
 \advance\countxviii@\m@ne\repeat                                           
 \dummyft@\tracinglostchars\z@\nopages@\frenchspacing\hbadness\@M}
\def\first@#1#2\end{#1}
\def\printoptions{\W@{Do you want S(yntax check),
  G(alleys) or P(ages)?}%
 \message{Type S, G or P, followed by <return>: }%
 \begingroup 
 \endlinechar\m@ne 
 \read\m@ne to\ans@
 \edef\ans@{\uppercase{\def\noexpand\ans@{%
   \expandafter\first@\ans@ P\end}}}%
 \expandafter\endgroup\ans@
 \if\ans@ P
 \else \if\ans@ S\syntax
 \else \if\ans@ G\galleys
 \else\message{? Unknown option: \ans@; using the `pages' option.}%
 \fi\fi\fi}
\def\alloc@#1#2#3#4#5{\global\advance\count1#1by\@ne
 \ch@ck#1#4#2\allocationnumber=\count1#1
 \global#3#5=\allocationnumber
 \ifalloc@\wlog{\string#5=\string#2\the\allocationnumber}\fi}
\def\document{\def\alloclist@{}\def\fontlist@{}}
\let\enddocument\bye

\let\proclaim\undefined
\let\footnote\undefined
\let\=\undefined
\let\>\undefined

\catcode`\@=\active
\message{... finished}

\expandafter\ifx\csname mathdefs.tex\endcsname\relax
  \expandafter\gdef\csname mathdefs.tex\endcsname{}
\else \message{Hey!  Apparently you were trying to
  \string\input{mathdefs.tex} twice.   This does not make sense.} 
\errmessage{Please edit your file (probably \jobname.tex) and remove
any duplicate ``\string\input'' lines}\endinput\fi




\catcode`\X=12\catcode`\@=11

\def\n@wcount{\alloc@0\count\countdef\insc@unt}
\def\n@wwrite{\alloc@7\write\chardef\sixt@@n}
\def\n@wread{\alloc@6\read\chardef\sixt@@n}
\def\r@s@t{\relax}\def\v@idline{\par}\def\@mputate#1/{#1}
\def\l@c@l#1X{\firstpart.#1}\def\gl@b@l#1X{#1}\def\t@d@l#1X{{}}

\def\crossrefs#1{\ifx\all#1\let\tr@ce=\all\else\def\tr@ce{#1,}\fi
   \n@wwrite\cit@tionsout\openout\cit@tionsout=\jobname.cit 
   \write\cit@tionsout{\tr@ce}\expandafter\setfl@gs\tr@ce,}
\def\setfl@gs#1,{\def\@{#1}\ifx\@\empty\let\next=\relax
   \else\let\next=\setfl@gs\expandafter\xdef
   \csname#1tr@cetrue\endcsname{}\fi\next}
\def\m@ketag#1#2{\expandafter\n@wcount\csname#2tagno\endcsname
     \csname#2tagno\endcsname=0\let\tail=\all\xdef\all{\tail#2,}
   \ifx#1\l@c@l\let\tail=\r@s@t\xdef\r@s@t{\csname#2tagno\endcsname=0\tail}\fi
   \expandafter\gdef\csname#2cite\endcsname##1{\expandafter
     \ifx\csname#2tag##1\endcsname\relax?\else\csname#2tag##1\endcsname\fi
     \expandafter\ifx\csname#2tr@cetrue\endcsname\relax\else
     \write\cit@tionsout{#2tag ##1 cited on page \folio.}\fi}
   \expandafter\gdef\csname#2page\endcsname##1{\expandafter
     \ifx\csname#2page##1\endcsname\relax?\else\csname#2page##1\endcsname\fi
     \expandafter\ifx\csname#2tr@cetrue\endcsname\relax\else
     \write\cit@tionsout{#2tag ##1 cited on page \folio.}\fi}
   \expandafter\gdef\csname#2tag\endcsname##1{\expandafter
      \ifx\csname#2check##1\endcsname\relax
      \expandafter\xdef\csname#2check##1\endcsname{}%
      \else\immediate\write16{Warning: #2tag ##1 used more than once.}\fi
      \multit@g{#1}{#2}##1/X%
      \write\t@gsout{#2tag ##1 assigned number \csname#2tag##1\endcsname\space
      on page \number\count0.}%
   \csname#2tag##1\endcsname}}

\def\multit@g#1#2#3/#4X{\def\t@mp{#4}\ifx\t@mp\empty%
      \global\advance\csname#2tagno\endcsname by 1 
      \expandafter\xdef\csname#2tag#3\endcsname
      {#1\number\csname#2tagno\endcsnameX}%
   \else\expandafter\ifx\csname#2last#3\endcsname\relax
      \expandafter\n@wcount\csname#2last#3\endcsname
      \global\advance\csname#2tagno\endcsname by 1 
      \expandafter\xdef\csname#2tag#3\endcsname
      {#1\number\csname#2tagno\endcsnameX}
      \write\t@gsout{#2tag #3 assigned number \csname#2tag#3\endcsname\space
      on page \number\count0.}\fi
   \global\advance\csname#2last#3\endcsname by 1
   \def\t@mp{\expandafter\xdef\csname#2tag#3/}%
   \expandafter\t@mp\@mputate#4\endcsname
   {\csname#2tag#3\endcsname\lastpart{\csname#2last#3\endcsname}}\fi}
\def\t@gs#1{\def\all{}\m@ketag#1e\m@ketag#1s\m@ketag\t@d@l p
\let\realscite\scite
\let\realstag\stag
   \m@ketag\gl@b@l r \n@wread\t@gsin
   \openin\t@gsin=\jobname.tgs \re@der \closein\t@gsin
   \n@wwrite\t@gsout\openout\t@gsout=\jobname.tgs }
\outer\def\localtags{\t@gs\l@c@l}
\outer\def\globaltags{\t@gs\gl@b@l}
\outer\def\newlocaltag#1{\m@ketag\l@c@l{#1}}
\outer\def\newglobaltag#1{\m@ketag\gl@b@l{#1}}

\newif\ifpr@ 
\def\m@kecs #1tag #2 assigned number #3 on page #4.%
   {\expandafter\gdef\csname#1tag#2\endcsname{#3}
   \expandafter\gdef\csname#1page#2\endcsname{#4}
   \ifpr@\expandafter\xdef\csname#1check#2\endcsname{}\fi}
\def\re@der{\ifeof\t@gsin\let\next=\relax\else
   \read\t@gsin to\t@gline\ifx\t@gline\v@idline\else
   \expandafter\m@kecs \t@gline\fi\let \next=\re@der\fi\next}
\def\pretags#1{\pr@true\pret@gs#1,,}
\def\pret@gs#1,{\def\@{#1}\ifx\@\empty\let\n@xtfile=\relax
   \else\let\n@xtfile=\pret@gs \openin\t@gsin=#1.tgs \message{#1} \re@der 
   \closein\t@gsin\fi \n@xtfile}

\newcount\sectno\sectno=0\newcount\subsectno\subsectno=0
\newif\ifultr@local \def\ultralocal{\ultr@localtrue}
\def\firstpart{\number\sectno}
\def\lastpart#1{\ifcase#1 \or a\or b\or c\or d\or e\or f\or g\or h\or 
   i\or k\or l\or m\or n\or o\or p\or q\or r\or s\or t\or u\or v\or w\or 
   x\or y\or z \fi}

\def\resetall{\global\advance\sectno by 1\subsectno=0
   \gdef\firstpart{\number\sectno}\r@s@t}
\def\resetsub{\global\advance\subsectno by 1
   \gdef\firstpart{\number\sectno.\number\subsectno}\r@s@t}
\def\newsection#1\par{\resetall\vskip0pt plus.3\vsize\penalty-250
   \vskip0pt plus-.3\vsize\bigskip\bigskip
   \message{#1}\leftline{\bf#1}\nobreak\bigskip}
\def\subsection#1\par{\ifultr@local\resetsub\fi
   \vskip0pt plus.2\vsize\penalty-250\vskip0pt plus-.2\vsize
   \bigskip\smallskip\message{#1}\leftline{\bf#1}\nobreak\medskip}


\newdimen\marginshift

\newdimen\margindelta
\newdimen\marginmax
\newdimen\marginmin

\def\margininit{       
\marginmax=3 true cm                  
				      
\margindelta=0.1 true cm              
\marginmin=0.1true cm                 
\marginshift=\marginmin
}    

\def\t@gsjj#1,{\def\@{#1}\ifx\@\empty\let\next=\relax\else\let\next=\t@gsjj
   \def\@@{p}\ifx\@\@@\else
   \expandafter\gdef\csname#1cite\endcsname##1{\citejj{##1}}
   \expandafter\gdef\csname#1page\endcsname##1{?}
   \expandafter\gdef\csname#1tag\endcsname##1{\tagjj{##1}}\fi\fi\next}
\newif\ifshowstuffinmargin
\showstuffinmarginfalse
\def\jjtags{\ifx\shlhetal\relax 
  \else
\ifx\shlhetal\undefinedcontrolseq
\else
\showstuffinmargintrue
\ifx\all\relax\else\expandafter\t@gsjj\all,\fi\fi \fi
}

\def\tagjj#1{\realstag{#1}\mginpar{\zeigen{#1}}}
\def\citejj#1{\rechnen{#1}\mginpar{\zeigen{#1}}}     

\def\rechnen#1{\expandafter\ifx\csname stag#1\endcsname\relax ??\else
                           \csname stag#1\endcsname\fi}

\newdimen\theight

\def\marginfont{\sevenrm}

\def\trymarginbox#1{\setbox0=\hbox{\marginfont\hskip\marginshift #1}%
		\global\marginshift\wd0 
		\global\advance\marginshift\margindelta}

\def \mginpar#1{%
\ifvmode\setbox0\hbox to \hsize{\hfill\rlap{\marginfont\quad#1}}%
\ht0 0cm
\dp0 0cm
\box0\vskip-\baselineskip
\else 
             \vadjust{\trymarginbox{#1}%
		\ifdim\marginshift>\marginmax \global\marginshift\marginmin
			\trymarginbox{#1}%
                \fi
             \theight=\ht0
             \advance\theight by \dp0    \advance\theight by \lineskip
             \kern -\theight \vbox to \theight{\rightline{\rlap{\box0}}%
\vss}}\fi}


\def\t@gsoff#1,{\def\@{#1}\ifx\@\empty\let\next=\relax\else\let\next=\t@gsoff
   \def\@@{p}\ifx\@\@@\else
   \expandafter\gdef\csname#1cite\endcsname##1{\zeigen{##1}}
   \expandafter\gdef\csname#1page\endcsname##1{?}
   \expandafter\gdef\csname#1tag\endcsname##1{\zeigen{##1}}\fi\fi\next}
\def\verbatimtags{\showstuffinmarginfalse
\ifx\all\relax\else\expandafter\t@gsoff\all,\fi}
\def\zeigen#1{\hbox{$\langle$}#1\hbox{$\rangle$}}

\def\margincite#1{\ifshowstuffinmargin\mginpar{\zeigen{#1}}\fi}

\def\(#1){\edef\dot@g{\ifmmode\ifinner(\hbox{\noexpand\etag{#1}})
   \else\noexpand\eqno(\hbox{\noexpand\etag{#1}})\fi
   \else(\noexpand\ecite{#1})\fi}\dot@g}

\newif\ifbr@ck
\def\eat#1{}
\def\[#1]{\br@cktrue[\br@cket#1'X]}
\def\br@cket#1'#2X{\def\temp{#2}\ifx\temp\empty\let\next\eat
   \else\let\next\br@cket\fi
   \ifbr@ck\br@ckfalse\br@ck@t#1,X\else\br@cktrue#1\fi\next#2X}
\def\br@ck@t#1,#2X{\def\temp{#2}\ifx\temp\empty\let\neext\eat
   \else\let\neext\br@ck@t\def\temp{,}\fi
   \def\teemp{#1}\ifx\teemp\empty\else\rcite{#1}\fi\temp\neext#2X}
\def\resetbr@cket{\gdef\[##1]{[\rtag{##1}]}}
\def\references{\resetbr@cket\newsection References\par}

\newtoks\symb@ls\newtoks\s@mb@ls\newtoks\p@gelist\n@wcount\ftn@mber
    \ftn@mber=1\newif\ifftn@mbers\ftn@mbersfalse\newif\ifbyp@ge\byp@gefalse
\def\defm@rk{\ifftn@mbers\n@mberm@rk\else\symb@lm@rk\fi}
\def\n@mberm@rk{\xdef\m@rk{{\the\ftn@mber}}%
    \global\advance\ftn@mber by 1 }
\def\rot@te#1{\let\temp=#1\global#1=\expandafter\r@t@te\the\temp,X}
\def\r@t@te#1,#2X{{#2#1}\xdef\m@rk{{#1}}}
\def\b@@st#1{{$^{#1}$}}\def\str@p#1{#1}
\def\symb@lm@rk{\ifbyp@ge\rot@te\p@gelist\ifnum\expandafter\str@p\m@rk=1 
    \s@mb@ls=\symb@ls\fi\write\f@nsout{\number\count0}\fi \rot@te\s@mb@ls}
\def\byp@ge{\byp@getrue\n@wwrite\f@nsin\openin\f@nsin=\jobname.fns 
    \n@wcount\currentp@ge\currentp@ge=0\p@gelist={0}
    \re@dfns\closein\f@nsin\rot@te\p@gelist
    \n@wread\f@nsout\openout\f@nsout=\jobname.fns }
\def\m@kelist#1X#2{{#1,#2}}
\def\re@dfns{\ifeof\f@nsin\let\next=\relax\else\read\f@nsin to \f@nline
    \ifx\f@nline\v@idline\else\let\t@mplist=\p@gelist
    \ifnum\currentp@ge=\f@nline
    \global\p@gelist=\expandafter\m@kelist\the\t@mplistX0
    \else\currentp@ge=\f@nline
    \global\p@gelist=\expandafter\m@kelist\the\t@mplistX1\fi\fi
    \let\next=\re@dfns\fi\next}
\def\symbols#1{\symb@ls={#1}\s@mb@ls=\symb@ls} 
\def\bigsymbol{\textstyle}
\symbols{\bigsymbol\ast,\dagger,\ddagger,\sharp,\flat,\natural,\star}
\def\ftnumbers{\ftn@mberstrue} \def\ftsymbols{\ftn@mbersfalse}
\def\paginal{\byp@ge} \def\resetftnumbers{\ftn@mber=1}
\def\ftnote#1{\defm@rk\expandafter\expandafter\expandafter\footnote
    \expandafter\b@@st\m@rk{#1}}

\long\def\jump#1\endjump{}
\def\ssum{\mathop{\lower .1em\hbox{$\textstyle\Sigma$}}\nolimits}

\def\qed{\nobreak\kern 1em \vrule height .5em width .5em depth 0em}
\def\newneq{\hbox{\rlap{\hbox to 1\wd9{\hss$=$\hss}}\raise .1em 
   \hbox to 1\wd9{\hss$\scriptscriptstyle/$\hss}}}
\def\subsetne{\setbox9 = \hbox{$\subset$}\mathrel{\hbox{\rlap
   {\lower .4em \newneq}\raise .13em \hbox{$\subset$}}}}
\def\supsetne{\setbox9 = \hbox{$\subset$}\mathrel{\hbox{\rlap
   {\lower .4em \newneq}\raise .13em \hbox{$\supset$}}}}

\def\vbar{\mathchoice{\vrule height6.3ptdepth-.5ptwidth.8pt\kern-.8pt}
   {\vrule height6.3ptdepth-.5ptwidth.8pt\kern-.8pt}
   {\vrule height4.1ptdepth-.35ptwidth.6pt\kern-.6pt}
   {\vrule height3.1ptdepth-.25ptwidth.5pt\kern-.5pt}}
\def\f@dge{\mathchoice{}{}{\mkern.5mu}{\mkern.8mu}}
\def\b@c#1#2{{\rm \mkern#2mu\vbar\mkern-#2mu#1}}
\def\b@b#1{{\rm I\mkern-3.5mu #1}}
\def\b@a#1#2{{\rm #1\mkern-#2mu\f@dge #1}}
\def\bb#1{{\count4=`#1 \advance\count4by-64 \ifcase\count4\or\b@a A{11.5}\or
   \b@b B\or\b@c C{5}\or\b@b D\or\b@b E\or\b@b F \or\b@c G{5}\or\b@b H\or
   \b@b I\or\b@c J{3}\or\b@b K\or\b@b L \or\b@b M\or\b@b N\or\b@c O{5} \or
   \b@b P\or\b@c Q{5}\or\b@b R\or\b@a S{8}\or\b@a T{10.5}\or\b@c U{5}\or
   \b@a V{12}\or\b@a W{16.5}\or\b@a X{11}\or\b@a Y{11.7}\or\b@a Z{7.5}\fi}}

\catcode`\X=11 \catcode`\@=12




\let\thischap\jobname

\def\partof#1{\csname returnthe#1part\endcsname}
\def\chapof#1{\csname returnthe#1chap\endcsname}

\def\setchapter#1,#2,#3.{%
  \expandafter\def\csname returnthe#1part\endcsname{#2}%
  \expandafter\def\csname returnthe#1chap\endcsname{#3}%
}

\setchapter 300a,A,I.
\setchapter 300b,A,II.
\setchapter 300c,A,III.
\setchapter 300d,A,IV.
\setchapter 300e,A,V.
\setchapter 300f,A,VI.
\setchapter 300g,A,VII.
\setchapter   88,B,I.
\setchapter  600,B,II.
\setchapter  705,B,III.

\def\cprefix#1{
\edef\theotherpart{\partof{#1}}\edef\theotherchap{\chapof{#1}}%
\ifx\theotherpart\thispart
   \ifx\theotherchap\thischap 
    \else 
     \theotherchap%
    \fi
   \else 
     \theotherpart.\theotherchap\fi}

\def\sectioncite[#1]#2{%
     \cprefix{#2}#1}

\edef\thispart{\partof{\thischap}}
\edef\thischap{\chapof{\thischap}}



\expandafter\ifx\csname citeadd.tex\endcsname\relax
\expandafter\gdef\csname citeadd.tex\endcsname{}
\else \message{Hey!  Apparently you were trying to
\string\input{citeadd.tex} twice.   This does not make sense.} 
\errmessage{Please edit your file (probably \jobname.tex) and remove
any duplicate ``\string\input'' lines}\endinput\fi

\sectno=-1   
\localtags
\jjtags
\NoBlackBoxes
\define\mr{\medskip\roster}
\define\sn{\smallskip\noindent}
\define\mn{\medskip\noindent}
\define\bn{\bigskip\noindent}
\define\ub{\underbar}
\define\wilog{\text{without loss of generality}}
\define\ermn{\endroster\medskip\noindent}

\define\dbcu{\dsize\bigcup}
\define \nl{\newline}
\newbox\noforkbox \newdimen\forklinewidth
\forklinewidth=0.3pt   
\setbox0\hbox{$\textstyle\bigcup$}
\setbox1\hbox to \wd0{\hfil\vrule width \forklinewidth depth \dp0
                        height \ht0 \hfil}
\wd1=0 cm
\setbox\noforkbox\hbox{\box1\box0\relax}
\def\unionstick{\mathop{\copy\noforkbox}\limits}
\def\nonfork#1#2_#3{#1\unionstick_{\textstyle #3}#2}
\def\nonforkin#1#2_#3^#4{#1\unionstick_{\textstyle #3}^{\textstyle #4}#2}     
%
\setbox0\hbox{$\textstyle\bigcup$}
\setbox1\hbox to \wd0{\hfil{\sl /\/}\hfil}
\setbox2\hbox to \wd0{\hfil\vrule height \ht0 depth \dp0 width
                                \forklinewidth\hfil}
\wd1=0cm
\wd2=0cm
\newbox\doesforkbox
\setbox\doesforkbox\hbox{\box1\box0\relax}
\def\nunionstick{\mathop{\copy\doesforkbox}\limits}

\def\fork#1#2_#3{#1\nunionstick_{\textstyle #3}#2}
\def\forkin#1#2_#3^#4{#1\nunionstick_{\textstyle #3}^{\textstyle #4}#2}     
\magnification=\magstep 1
\documentstyle{amsppt}

{    
\catcode`@11

\ifx\alicetwothousandloaded@\relax
  \endinput\else\global\let\alicetwothousandloaded@\relax\fi

\gdef\subjclass{\let\savedef@\subjclass
 \def\subjclass##1\endsubjclass{\let\subjclass\savedef@
   \toks@{\def\usualspace{{\rm\enspace}}\eightpoint}%
   \toks@@{##1\unskip.}%
   \edef\thesubjclass@{\the\toks@
     \frills@{{\noexpand\rm2000 {\noexpand\it Mathematics Subject
       Classification}.\noexpand\enspace}}%
     \the\toks@@}}%
  \nofrillscheck\subjclass}
} 


\expandafter\ifx\csname alice2jlem.tex\endcsname\relax
  \expandafter\xdef\csname alice2jlem.tex\endcsname{\the\catcode`@}
\else \message{Hey!  Apparently you were trying to
\string\input{alice2jlem.tex}  twice.   This does not make sense.}
\errmessage{Please edit your file (probably \jobname.tex) and remove
any duplicate ``\string\input'' lines}\endinput\fi

\expandafter\ifx\csname bib4plain.tex\endcsname\relax
  \expandafter\gdef\csname bib4plain.tex\endcsname{}
\else \message{Hey!  Apparently you were trying to \string\input
  bib4plain.tex twice.   This does not make sense.}
\errmessage{Please edit your file (probably \jobname.tex) and remove
any duplicate ``\string\input'' lines}\endinput\fi

\def\renewcommand{\newcommand}	       
\edef\cite{\the\catcode`@}%
\catcode`@ = 11
\let\@oldatcatcode = \cite
\chardef\@letter = 11
\chardef\@other = 12
%
%
%
%
\def\@innerdef#1#2{\edef#1{\expandafter\noexpand\csname #2\endcsname}}%
%
%
\@innerdef\@innernewcount{newcount}%
\@innerdef\@innernewdimen{newdimen}%
\@innerdef\@innernewif{newif}%
\@innerdef\@innernewwrite{newwrite}%
%
%
%
\def\@gobble#1{}%
%
%
%
\ifx\inputlineno\@undefined
   \let\@linenumber = \empty 
\else
   \def\@linenumber{\the\inputlineno:\space}%
\fi
%
%
%
\def\@futurenonspacelet#1{\def\cs{#1}%
   \afterassignment\@stepone\let\@nexttoken=
}%
\begingroup 
\def\\{\global\let\@stoken= }%
\\ 
\endgroup
\def\@stepone{\expandafter\futurelet\cs\@steptwo}%
\def\@steptwo{\expandafter\ifx\cs\@stoken\let\@@next=\@stepthree
   \else\let\@@next=\@nexttoken\fi \@@next}%
\def\@stepthree{\afterassignment\@stepone\let\@@next= }%
%
%
%
\def\@getoptionalarg#1{%
   \let\@optionaltemp = #1%
   \let\@optionalnext = \relax
   \@futurenonspacelet\@optionalnext\@bracketcheck
}%
%
%
\def\@bracketcheck{%
   \ifx [\@optionalnext
      \expandafter\@@getoptionalarg
   \else
      \let\@optionalarg = \empty
      \expandafter\@optionaltemp
   \fi
}%
\def\@@getoptionalarg[#1]{%
   \def\@optionalarg{#1}%
   \@optionaltemp
}%
%
%
%
\def\@nnil{\@nil}%
\def\@fornoop#1\@@#2#3{}%
\def\@for#1:=#2\do#3{%
   \edef\@fortmp{#2}%
   \ifx\@fortmp\empty \else
      \expandafter\@forloop#2,\@nil,\@nil\@@#1{#3}%
   \fi
}%
\def\@forloop#1,#2,#3\@@#4#5{\def#4{#1}\ifx #4\@nnil \else
       #5\def#4{#2}\ifx #4\@nnil \else#5\@iforloop #3\@@#4{#5}\fi\fi
}%
\def\@iforloop#1,#2\@@#3#4{\def#3{#1}\ifx #3\@nnil
       \let\@nextwhile=\@fornoop \else
      #4\relax\let\@nextwhile=\@iforloop\fi\@nextwhile#2\@@#3{#4}%
}%
%
%
%
\@innernewif\if@fileexists
\def\@testfileexistence{\@getoptionalarg\@finishtestfileexistence}%
\def\@finishtestfileexistence#1{%
   \begingroup
      \def\extension{#1}%
      \immediate\openin0 =
         \ifx\@optionalarg\empty\jobname\else\@optionalarg\fi
         \ifx\extension\empty \else .#1\fi
         \space
      \ifeof 0
         \global\@fileexistsfalse
      \else
         \global\@fileexiststrue
      \fi
      \immediate\closein0
   \endgroup
}%
%
%
%
%
\def\bibliographystyle#1{%
   \@readauxfile
   \@writeaux{\string\bibstyle{#1}}%
}%
\let\bibstyle = \@gobble
%
%
\let\bblfilebasename = \jobname
\def\bibliography#1{%
   \@readauxfile
   \@writeaux{\string\bibdata{#1}}%
   \@testfileexistence[\bblfilebasename]{bbl}%
   \if@fileexists
      \nobreak
      \@readbblfile
   \fi
}%
\let\bibdata = \@gobble
%
%
\def\nocite#1{%
   \@readauxfile
   \@writeaux{\string\citation{#1}}%
}%
\@innernewif\if@notfirstcitation
%
%
\def\cite{\@getoptionalarg\@cite}%
%
%
\def\@cite#1{%
   \let\@citenotetext = \@optionalarg
   \printcitestart
   \nocite{#1}%
   \@notfirstcitationfalse
   \@for \@citation :=#1\do
   {%
      \expandafter\@onecitation\@citation\@@
   }%
   \ifx\empty\@citenotetext\else
      \printcitenote{\@citenotetext}%
   \fi
   \printcitefinish
}%
\newif\ifweareinprivate
\weareinprivatetrue
\ifx\shlhetal\undefinedcontrolseq\weareinprivatefalse\fi
\ifx\shlhetal\relax\weareinprivatefalse\fi
\def\@onecitation#1\@@{%
   \if@notfirstcitation
      \printbetweencitations
   \fi
   \expandafter \ifx \csname\@citelabel{#1}\endcsname \relax
      \if@citewarning
         \message{\@linenumber Undefined citation `#1'.}%
      \fi
     \ifweareinprivate
      \expandafter\gdef\csname\@citelabel{#1}\endcsname{%
\strut 
\vadjust{\vskip-\dp\strutbox
\vbox to 0pt{\vss\parindent0cm \leftskip=\hsize 
\advance\leftskip3mm
\advance\hsize 4cm\strut\openup-4pt 
\rightskip 0cm plus 1cm minus 0.5cm ?  #1 ?\strut}}
         {\tt
            \escapechar = -1
            \nobreak\hskip0pt\pfeilsw
            \expandafter\string\csname#1\endcsname
             \pfeilso
            \nobreak\hskip0pt
         }%
      }%
     \else  
      \expandafter\gdef\csname\@citelabel{#1}\endcsname{%
            {\tt\expandafter\string\csname#1\endcsname}
      }%
     \fi  
   \fi
   \csname\@citelabel{#1}\endcsname
   \@notfirstcitationtrue
}%
%
%
\def\@citelabel#1{b@#1}%
%
%
\def\@citedef#1#2{\expandafter\gdef\csname\@citelabel{#1}\endcsname{#2}}%
%
%
%
\def\@readbblfile{%
   \ifx\@itemnum\@undefined
      \@innernewcount\@itemnum
   \fi
   \begingroup
      \def\begin##1##2{%
         \setbox0 = \hbox{\biblabelcontents{##2}}%
         \biblabelwidth = \wd0
      }%
      \def\end##1{}
      %
      %
      \@itemnum = 0
      \def\bibitem{\@getoptionalarg\@bibitem}%
      \def\@bibitem{%
         \ifx\@optionalarg\empty
            \expandafter\@numberedbibitem
         \else
            \expandafter\@alphabibitem
         \fi
      }%
      \def\@alphabibitem##1{%
         \expandafter \xdef\csname\@citelabel{##1}\endcsname {\@optionalarg}%
         \ifx\biblabelprecontents\@undefined
            \let\biblabelprecontents = \relax
         \fi
         \ifx\biblabelpostcontents\@undefined
            \let\biblabelpostcontents = \hss
         \fi
         \@finishbibitem{##1}%
      }%
      \def\@numberedbibitem##1{%
         \advance\@itemnum by 1
         \expandafter \xdef\csname\@citelabel{##1}\endcsname{\number\@itemnum}%
         \ifx\biblabelprecontents\@undefined
            \let\biblabelprecontents = \hss
         \fi
         \ifx\biblabelpostcontents\@undefined
            \let\biblabelpostcontents = \relax
         \fi
         \@finishbibitem{##1}%
      }%
      \def\@finishbibitem##1{%
         \biblabelprint{\csname\@citelabel{##1}\endcsname}%
         \@writeaux{\string\@citedef{##1}{\csname\@citelabel{##1}\endcsname}}%
         \ignorespaces
      }%
      %
      %
      \let\em = \bblem
      \let\newblock = \bblnewblock
      \let\sc = \bblsc
      \frenchspacing
      \clubpenalty = 4000 \widowpenalty = 4000
      \tolerance = 10000 \hfuzz = .5pt
      \everypar = {\hangindent = \biblabelwidth
                      \advance\hangindent by \biblabelextraspace}%
      \bblrm
      \parskip = 1.5ex plus .5ex minus .5ex
      \biblabelextraspace = .5em
      \bblhook
      \input \bblfilebasename.bbl
   \endgroup
}%
%
%
\@innernewdimen\biblabelwidth
\@innernewdimen\biblabelextraspace
%
%
%
\def\biblabelprint#1{%
   \noindent
   \hbox to \biblabelwidth{%
      \biblabelprecontents
      \biblabelcontents{#1}%
      \biblabelpostcontents
   }%
   \kern\biblabelextraspace
}%
%
%
%
\def\biblabelcontents#1{{\bblrm [#1]}}%
%
%
\def\bblrm{\rm}%
%
%
\def\bblem{\it}%
%
%
\def\bblsc{\ifx\@scfont\@undefined
              \font\@scfont = cmcsc10
           \fi
           \@scfont
}%
%
%
\def\bblnewblock{\hskip .11em plus .33em minus .07em }%
%
%
\let\bblhook = \empty
%
%
%
\def\printcitestart{[}
\def\printcitefinish{]}
\def\printbetweencitations{, }
\def\printcitenote#1{, #1}
%
%
%
\let\citation = \@gobble
%
%
%
\@innernewcount\@numparams
%
%
\def\newcommand#1{%
   \def\@commandname{#1}%
   \@getoptionalarg\@continuenewcommand
}%
%
%
\def\@continuenewcommand{%
   \@numparams = \ifx\@optionalarg\empty 0\else\@optionalarg \fi \relax
   \@newcommand
}%
%
%
\def\@newcommand#1{%
   \def\@startdef{\expandafter\edef\@commandname}%
   \ifnum\@numparams=0
      \let\@paramdef = \empty
   \else
      \ifnum\@numparams>9
         \errmessage{\the\@numparams\space is too many parameters}%
      \else
         \ifnum\@numparams<0
            \errmessage{\the\@numparams\space is too few parameters}%
         \else
            \edef\@paramdef{%
               \ifcase\@numparams
                  \empty  No arguments.
               \or ####1%
               \or ####1####2%
               \or ####1####2####3%
               \or ####1####2####3####4%
               \or ####1####2####3####4####5%
               \or ####1####2####3####4####5####6%
               \or ####1####2####3####4####5####6####7%
               \or ####1####2####3####4####5####6####7####8%
               \or ####1####2####3####4####5####6####7####8####9%
               \fi
            }%
         \fi
      \fi
   \fi
   \expandafter\@startdef\@paramdef{#1}%
}%
%
%
%
%
\def\@readauxfile{%
   \if@auxfiledone \else 
      \global\@auxfiledonetrue
      \@testfileexistence{aux}%
      \if@fileexists
         \begingroup
            \endlinechar = -1
            \catcode`@ = 11
            \input \jobname.aux
         \endgroup
      \else
         \message{\@undefinedmessage}%
         \global\@citewarningfalse
      \fi
      \immediate\openout\@auxfile = \jobname.aux
   \fi
}%
%
%
\newif\if@auxfiledone
\ifx\noauxfile\@undefined \else \@auxfiledonetrue\fi
%
%
%
%
\@innernewwrite\@auxfile
\def\@writeaux#1{\ifx\noauxfile\@undefined \write\@auxfile{#1}\fi}%
%
%
%
\ifx\@undefinedmessage\@undefined
   \def\@undefinedmessage{No .aux file; I won't give you warnings about
                          undefined citations.}%
\fi
%
%
\@innernewif\if@citewarning
\ifx\noauxfile\@undefined \@citewarningtrue\fi
%
%
%
\catcode`@ = \@oldatcatcode

\def\pfeilso{\leavevmode
            \vrule width 1pt height9pt depth 0pt\relax
           \vrule width 1pt height8.7pt depth 0pt\relax
           \vrule width 1pt height8.3pt depth 0pt\relax
           \vrule width 1pt height8.0pt depth 0pt\relax
           \vrule width 1pt height7.7pt depth 0pt\relax
            \vrule width 1pt height7.3pt depth 0pt\relax
            \vrule width 1pt height7.0pt depth 0pt\relax
            \vrule width 1pt height6.7pt depth 0pt\relax
            \vrule width 1pt height6.3pt depth 0pt\relax
            \vrule width 1pt height6.0pt depth 0pt\relax
            \vrule width 1pt height5.7pt depth 0pt\relax
            \vrule width 1pt height5.3pt depth 0pt\relax
            \vrule width 1pt height5.0pt depth 0pt\relax
            \vrule width 1pt height4.7pt depth 0pt\relax
            \vrule width 1pt height4.3pt depth 0pt\relax
            \vrule width 1pt height4.0pt depth 0pt\relax
            \vrule width 1pt height3.7pt depth 0pt\relax
            \vrule width 1pt height3.3pt depth 0pt\relax
            \vrule width 1pt height3.0pt depth 0pt\relax
            \vrule width 1pt height2.7pt depth 0pt\relax
            \vrule width 1pt height2.3pt depth 0pt\relax
            \vrule width 1pt height2.0pt depth 0pt\relax
            \vrule width 1pt height1.7pt depth 0pt\relax
            \vrule width 1pt height1.3pt depth 0pt\relax
            \vrule width 1pt height1.0pt depth 0pt\relax
            \vrule width 1pt height0.7pt depth 0pt\relax
            \vrule width 1pt height0.3pt depth 0pt\relax}

\def\pfeilsw{ \leavevmode 
            \vrule width 1pt height0.3pt depth 0pt\relax
            \vrule width 1pt height0.7pt depth 0pt\relax
            \vrule width 1pt height1.0pt depth 0pt\relax
            \vrule width 1pt height1.3pt depth 0pt\relax
            \vrule width 1pt height1.7pt depth 0pt\relax
            \vrule width 1pt height2.0pt depth 0pt\relax
            \vrule width 1pt height2.3pt depth 0pt\relax
            \vrule width 1pt height2.7pt depth 0pt\relax
            \vrule width 1pt height3.0pt depth 0pt\relax
            \vrule width 1pt height3.3pt depth 0pt\relax
            \vrule width 1pt height3.7pt depth 0pt\relax
            \vrule width 1pt height4.0pt depth 0pt\relax
            \vrule width 1pt height4.3pt depth 0pt\relax
            \vrule width 1pt height4.7pt depth 0pt\relax
            \vrule width 1pt height5.0pt depth 0pt\relax
            \vrule width 1pt height5.3pt depth 0pt\relax
            \vrule width 1pt height5.7pt depth 0pt\relax
            \vrule width 1pt height6.0pt depth 0pt\relax
            \vrule width 1pt height6.3pt depth 0pt\relax
            \vrule width 1pt height6.7pt depth 0pt\relax
            \vrule width 1pt height7.0pt depth 0pt\relax
            \vrule width 1pt height7.3pt depth 0pt\relax
            \vrule width 1pt height7.7pt depth 0pt\relax
            \vrule width 1pt height8.0pt depth 0pt\relax
            \vrule width 1pt height8.3pt depth 0pt\relax
            \vrule width 1pt height8.7pt depth 0pt\relax
            \vrule width 1pt height9pt depth 0pt\relax
      }


\def\widestnumber#1#2{}

\def\citewarning#1{\ifx\shlhetal\relax 
    \else
    \par{#1}\par
    \fi
}

\def\rm{\fam0 \tenrm}

\def\fakesubhead#1\endsubhead{\bigskip\noindent{\bf#1}\par}



%
%
%

%

\font\textrsfs=rsfs10
\font\scriptrsfs=rsfs7
\font\scriptscriptrsfs=rsfs5

\newfam\rsfsfam
\textfont\rsfsfam=\textrsfs
\scriptfont\rsfsfam=\scriptrsfs
\scriptscriptfont\rsfsfam=\scriptscriptrsfs

\edef\oldcatcodeofat{\the\catcode`\@}
\catcode`\@11

\def\Cal@@#1{\noaccents@ \fam \rsfsfam #1}

\catcode`\@\oldcatcodeofat


\expandafter\ifx \csname margininit\endcsname \relax\else\margininit\fi

\long\def\red#1\endred{}
\long\def\green#1\endgreen{}
\long\def\blue#1\endblue{}

\def\endred{ \unmatched endred! }
\def\endgreen{ \unmatched endgreen! }
\def\endblue{ \unmatched endblue! }

\ifx\latexcolors\undefinedcs\def\latexcolors{}\fi

\def\emptycs{}
\def\evaluatelatexcolors{%
        \ifx\latexcolors\emptycs\else
        \expandafter\xxevaluate\latexcolors\xxfertig\evaluatelatexcolors\fi}
\def\xxevaluate#1,#2\xxfertig{\setupthiscolor{#1}%
        \def\latexcolors{#2}}

\font\smallfont=cmsl7
\def\rutgerscolor{\ifmmode\else\endgraf\fi\smallfont
\advance\leftskip0.5cm\relax}
\def\setupthiscolor#1{\edef\tmptmpcs{\noexpand\bgroup\noexpand\rutgerscolor
\noexpand\def\noexpand\currentcolor{#1}%
\noexpand}%
\expandafter\let\csname#1\endcsname\tmptmpcs
\def\tmptmpcs{\checkColorUnmatched{#1}\popthecolor}
\expandafter\let\csname end#1\endcsname\tmptmpcs}

\def\checkColorUnmatched#1{\def\expectcolor{#1}%
    \ifx\expectcolor\currentcolor   
    \else \edef\failhere{\noexpand\tryingToClose '\currentcolor' with end\expectcolor}\failhere\fi}

\def\currentcolor{???}

\def\popthecolor{\ifmmode\else\endgraf\fi\egroup}

\expandafter\def\csname#1\endcsname{}

\evaluatelatexcolors

 \let\outerhead\head
 \def\head{\innerhead}
 \let\innerhead\outerhead

 \let\outersubhead\subhead
 \def\subhead{\innersubhead}
 \let\innersubhead\outersubhead

 \let\outersubsubhead\subsubhead
 \def\subsubhead{\innersubsubhead}
 \let\innersubsubhead\outersubsubhead

 \def\proclaim{\innerproclaim}
 \let\innerproclaim\outerproclaim

 %
 %
 %
 %

\def\demo#1{\medskip\noindent{\it #1.\/}}
\def\enddemo{\smallskip}

\def\remark#1{\medskip\noindent{\it #1.\/}}
\def\endremark{\smallskip}

\pageheight{8.5truein}
\topmatter
\title\nofrills{ON MODEL COMPLETION OF $T_{\text{aut}}$} \endtitle
\author {Saharon Shelah \thanks {\null\newline 
Partially supported by the Israel Science Foundation. Publication E34. \null\newline 
I would like to thank 
Alice Leonhardt for the beautiful typing. } \endthanks} \endauthor 

\affil{Institute of Mathematics\\
 The Hebrew University\\
 Jerusalem, Israel
 \medskip
 Rutgers University\\
 Mathematics Department\\
 New Brunswick, NJ  USA} \endaffil
\endtopmatter
\document  
 
\newpage

\head {Annotated Content} \endhead  \resetall 
\bn
\S0 Introduction
\mn
\S1
\mr
\item "{${{}}$}"  [We characterize stable $T$ for which the model
completion of $T_{\text{aut}}$ is stable (i.e., every
completion is).]
\ermn
\S2
\mr
\item "{${{}}$}"  [We prove that ``some completion is stable" is
  different and characterize it.]
\ermn
\S3
\mr
\item "{${{}}$}"  [We prove that if $T$ is stable, $T_{\text{aut}}$
has a model completion, $T_*$ is an unstable complete of
$T^{\text{mc}}_{\text{aut}}$, \ub{then} $T_*$ satisfies NSOP$_3$.
Moreover, simplicity is preserved.]
\endroster
\newpage

\head {\S0 Introduction} \endhead  \resetall \sectno=0
\bigskip

On the subject, history and background see \cite{BlSh:759}. 
For a complete first order $T$ they dealt with the existence of the
model completion $T_{\text{aut}}$ of 
$T \cup \{\sigma$ is an automorphism (for
$\tau_T$)$\}$. \nl
We may ask:
\nl
\ub{\stag{0.0} Question}:  If $T$ is stable and $T_{\text{aut}}$
has model completion $T^{\text{mc}}_{\text{aut}}$, when is (every)
completion of $T^{\text{mc}}_{\text{aut}}$ stable?

We answer in \scite{1.1} (observation \scite{1.2} deals with some
obvious things).
\nl
Section 1 raises some question which we discuss below (assuming $T$
stable, $T^{\text{mc}}_{\text{aut}}$ exists) some of which are
answered below.
\nl
\ub{\stag{0.1} Question}:  1) Can we in Claim \scite{1.1} below replace ``every
completion of $T^{\text{mc}}_{\text{aut}}$ is stable" by ``some completion
of $T^{\text{mc}}_{\text{aut}}$ is stable"? \nl
2) The ``unstable" in \scite{1.1} clause (a) can be replaced by
``having the independence property"; but can $T^{\text{mc}}_{\text{aut}}$ be
completed to a theory with the strict order property?  The
SOP$_n$'s?
\nl
3) What occurs if $T^{\text{mc}}_{\text{aut}}$ does not exist, can we still
say something? \nl
4) Point out that $(a)(\equiv(b))$ of \scite{1.1} holds (for some
stable $T$ for which $T^{\text{mc}}_{\text{aut}}$ exists) and fails for
others. \nl
5) Show for stable $T$ with $T^{\text{mc}}_{\text{aut}}$, that no completion
$T_*$ of $T^{\text{dc}}_{\text{ut}}$ has the explicit ncp (which means
that for some first order $E(\bar x,\bar y,\bar z)$, for every $n$ for
some $\bar c \subseteq {\frak C},
E(\bar x,\bar y,\bar c)$ is an equivalent relation which 
has $\ge n, <\aleph_0$ equivalence classes); \nl
a stronger version is \nl
6) For such $T,T_*$ can $T_*$ have obstructions (see \S4)? \nl
7) What if we use $\sigma_1,\sigma_2$?  
What about $\sigma_1,\dotsc,\sigma_n$?  What about pairwise commuting
$\sigma_1,\dotsc,\sigma_n$?  This is like
$(T_{\text{aut}})_{\text{aut}}$ for $n=2$.  \nl
8) Is there unstable $T$ such that $T_{\text{aut}}$ has model
completion?  (A conjecture stating that had been the starting point of Kikyo
Shelah \cite{KkSh:748}).
\bn
\ub{\stag{0.2} Discussion}:  We prove that:
\mr
\item "{$(A)$}"  on \scite{0.1}(1), for some $T$ (stable with
$T^{\text{mc}}_{\text{aut}}$ existing), some completion of
$T^{\text{mc}}_{\text{aut}}$ are stable and some are not (still we may
wonder on a general characerization, see \scite{2.4} below)
\sn
\item "{$(B)$}"  we shall show that for no such $T$ is any completion
of $T^{\text{mc}}_{\text{aut}}$ with the strict order property and even have
NSOP$_3$, see \scite{3.1}
\sn
\item "{$(C)$}"  we can look at the class of existentially closed
models of $T_{\text{aut}}$ (see \cite{ShUs:789} and
references there); the results are similar.
\endroster
\bn
Note
\demo{\stag{0.3} Observation} [Here?]
\mr
\item "{$(\alpha)$}"  for $T =$ theory of equality, $T_{\text{aut}}$
has a model completion and all completion
of $T^{\text{mc}}_{\text{aut}}$ are stable
\sn
\item "{$(\beta)$}"  for $T$ from \scite{2.2}, some completions of
$T^{\text{mc}}_{\text{aut}}$ are stable and some are not
\sn
\item "{$(\gamma)$}"  for $T = \text{ Th}(M \restriction
\{E,F_1,F_2,Q\}),M$ from \scite{2.2}, we get that all the completions
of $T^{\text{mc}}_{\text{aut}}$ are unstable.
\endroster
\enddemo
\bn
I think \nl
\ub{\stag{0.4} Quesiton}:  What about getting (in \S3) that
\mr
\item "{$(a)$}"  $T^{\text{mc}}_{\text{aut}}$ is simple in \S3?
\sn
\item "{$(b)$}"  even if $T$ is just simple,
$T^{\text{mc}}_{\text{aut}} \models$ NSOP$_3$
\sn
\item "{$(c)$}"  non elementary class (true).
\ermn
See below.
\newpage

\head {\S1 On the stability of 
model completion for $T_{\text{aut}}$ ($=T+ \sigma$ an automorphism)} \endhead  \resetall 
\bigskip

\demo{\stag{1.0} Hypothesis}  1) $T$ is first order complete and for
notational simplicity every formula is equivalent to a relation and
$\tau_T$ having only predicates. \nl
2) ${\frak C}$ is the monster model of $T$.
\enddemo
\bigskip

\definition{\stag{1.0a} Definition}  1) $T_{\text{aut}}$ is 
$T \cup \{\sigma$ is an automorphism (for $\tau_T)\}$, so $\sigma$ 
is a new unary function symbol that is $T_{\text{aut}} = T \cup
\{(\forall x_0,\dotsc,x_{n-1})[R(x_0,\dotsc,x_{n-1}) \equiv
R(\sigma(x_0),\dotsc,\sigma(x_{n-1}))]:R$ an $n$-place predicate of
$\tau_T\}$.
\nl
2) $T^{\text{mc}}_{\text{aut}}$ is the model completion, if it
exists. \nl
3) Let $T_*$ denote any completion of $T^{\text{mc}}_{\text{aut}}$ 
and $\sigma_*$ or
$\sigma^{N^+}$ is an automorphism. \nl
4) A completion $T_*$ of $T^{\text{mc}}_{\text{aut}}$ is cute \ub{if} it
has a model $N^+$ such that for some $M^+ \subseteq N^+$ we have
$\sigma^{N^+} = \text{ id}_{N^+}$.
\enddefinition
\bigskip

\definition{\stag{1.1b} Definition}  For $T$ as in \scite{0.1} let:
\nl
1) $K_{\text{aut}}(T) =$ the class of models of $T_{\text{aut}}$. \nl
2) $K^{\text{ec}}_{\text{aut}}(T) =$ the class of e.c. models of
$T_{\text{aut}}$. \nl
3) $K_*(T)$ is a subclass of $K^{\text{ec}}_{\text{aut}}(T)$ such that $M
\cong N \in K_* \Rightarrow M \in K_*$ and if $M \subseteq N$ are from
$K^{\text{ec}}_{\text{aut}}$ then $M \in K_* \Leftrightarrow N \in
K_*$; there are $\le 2^{|T|}$ such classes. \nl
4) $K_*$ is cute, etc. \nl
5) ${\frak C}_{\text{aut}}$ is a monster model for
$K^{\text{ec}}_{\text{aut}}$, i.e., a member of
$K^{\text{ec}}_{\text{aut}}$ which is $\bar \kappa$-saturated of
cardinality $\bar \kappa$; it is unique if $K_{\text{aut}}(T)$ has the JEP.
\nl
6) A class $K_*$ is stable \footnote{this is for classes as above, for
general non first order classes this does not fit}
\ub{if} for some $\lambda < \bar \kappa$ there is no
model $M \in K_*,m < \omega,\bar a_i \in {}^m M,i < \lambda$ and
q.f. formula $\varphi(\bar x,\bar y)$ which order $\{\bar a_i:i <
\lambda\}$. \nl
7) $K_*$ is simple \ub{if} there is a q.f. formula 
$\varphi(\bar x,\bar y)$ and $m$ such that for every $\lambda,\kappa$
we can find $M \in K_*,\bar a_\eta \i {}^{\ell g(\bar y)} M$ for $\eta
\in {}^{\kappa >} \lambda$ and $\bar b_\nu \in {}^{\ell g(\bar x)} M$
for $\nu \in {}^\kappa \lambda$ such that:
\mr
\item "{$(i)$}"  $M \models \varphi(\bar b_\eta,\bar a_{\eta
\restriction \alpha})$ for $\alpha < \kappa,\eta \in {}^\kappa
\lambda$
\sn
\item "{$(ii)$}"  no sequence in $m$ realizes $\ge m$ of the formulas
$\langle \varphi(\bar x,\bar a)_{\eta \char 94 <1>}:i < \lambda\}$.
\endroster
\enddefinition
\bn
On such models see \cite{Sh:54}, [xx], [xx].
\nl
\ub{\stag{1.1c} Fact}:  If $T^{\text{mc}}_{\text{aut}}$ exists then
$K^{\text{ec}}_{\text{aut}}(T)$ is the class of its models.
\bigskip

\proclaim{\stag{1.1y} Claim}  In the claims below we can replace ``$T$
has model completion" by dealing with the class
$K^{\text{ec}}_{\text{aut}}(T)$, and replace $T^*$ is a model
completion by dealing with $K_*$.
\endproclaim
\bigskip

\proclaim{\stag{1.1} Claim}  Let $T$ be stable, $T^{\text{mc}}_{\text{aut}}$
exists.  The $(a) \Leftrightarrow (b)$ where
\mr
\item "{$(a)$}"  $T^{\text{mc}}_{\text{aut}}$ is stable (i.e., every
completion is stable) 
\sn
\item "{$(b)$}"  if $M_0 \prec M_\ell \prec {\frak C}$ for
$\ell=1,2$ and $\nonfork{M_1}{M_2}_{M_0}$ then in ${\frak
C}^{\text{eq}},{\text{\rm acl\/}}_{{\frak C}^{\text{eq}}}
(M_1 \cup M_2) = { \text{\rm dcl\/}}_{{\frak C}^{\text{eq}}}(M_1 \cup
M_2)$
\sn
\item "{$(c)$}"  $T^{\text{mc}}_{\text{aut}}$ is dependent (i.e.,
every completion does not have the independence property).
\endroster
\endproclaim
\bigskip

\demo{Proof}  \ub{$(b) \Rightarrow (a)$}

We work in ${\frak C}^{\text{eq}}$ and use observation \scite{1.2} below.
Suppose ${\frak C}_* = ({\frak C},\sigma_*)$ is an expansion of
${\frak C}^{\text{eq}}$ to a model of
$T^{\text{mc}}_{\text{aut}}$ and let $\sigma^{\text{eq}}_*$
be the unique extension of $\sigma_*$ to an automorphism of ${\frak
C}^{\text{eq}}$. 
Let
$\lambda = \lambda^{|T|},M^+ \prec ({\frak C}^{\text{eq}},
\sigma^{\text{eq}}_*),|M^+|
= \lambda$ (note $|T| \ge \aleph_0$ here (by \scite{1.0}(1))).

Now for every $p \in {\Cal S}(M^+,{\frak C}^{\text{eq}}_*)$ let $a_p \in {\frak
C}$ realize $p$ in $({\frak C},\sigma_*)$ and let $M^+_p,N^+_p$ be such that

$$
M^+_p \prec M^+,\|M^+_p\| = |T| + \aleph_0
$$

$$
M^+_p \prec N^+_p \prec {\frak C}_\sigma,\|N^+_p\| = |T|
$$

$$
a_p \in N^+_p
$$

$$
\nonforkin{N^+_p \restriction \tau_T}{M^+ \restriction \tau_T}_{M^+_p
\restriction \tau_T}^{{\frak C}}.
$$
\mn
Let $A_p = \text{ acl}_{{\frak C}^{\text{eq}}}
(|N^+_p| \cup |M^+_p|)$.  We define a
two-place relation $E$ on ${\Cal S}(M^+,{\frak C}_\sigma)$ as follows:
\mr
\item "{$\circledast$}"  $p E q$
\ub{iff} $M^+_p = M^+_q$ and there is an isomorphism $f$ from $N^+_p$
onto $N^+_q$ which is the identity on $M^+_p$ and satisfying
$f_p(a_p) = a_q$.
\ermn
Clearly
\mr
\item "{$\circledast_0$}"  $E$ is an equivalence relation on 
${\Cal S}(M^+,{\frak C}_*)$
\sn
\item "{$\circledast_1$}"  $|{\Cal S}(M^+,{\frak C}_*)/E| \le
\lambda^{|T|}$.
\ermn
Hence it is enough to prove that
\mr
\item "{$\circledast_2$}"  $pEq \Rightarrow p=q$.
\endroster
\enddemo
\bigskip

\demo{Proof of $\circledast_2$}  Let $f$ witness $p Eq$, 

Let $f^+:A_p = \text{ dcl}_{{\frak C}^{\text{eq}}}(|N^+_p| \cup |M^+|)
\rightarrow A_q$ extends $f \cup \text{ id}_M$ 
and be an elementary mapping (in ${\frak C}^{\text{eq}}$); by non
forking calculus it exists and is unique.  Obviously
it commutes with $\sigma_*$.  Also $A_p$ (and $A_q$)
are algebraically closed sets in ${\frak C}^{\text{eq}}$ 
by our hypothesis (that is, clause (b)) applied to $|M^+_p|,|N^+_p|,|M^+|$
hence by \scite{1.2}(4), \scite{1.2t}(4) below, 
$f^+$ can be extended to an automorphism of ${\frak C}^{\text{eq}}$.
So by properties of model completion (and the obvious \scite{1.2t}(1) below)
we are done. 
\mn
\ub{$\neg(b) \Rightarrow \neg(a)$}:

Let $M_0,M_1,M_2$ form a counterexample to $(b)$.  
So let $e \in {\text{\rm acl\/}}_{{\frak C}^{\text{eq}}}
(M_1 \cup M_2) \backslash
{\text{\rm dcl\/}}_{{\frak C}^{\text{eq}}}(M_1 \cup M_2)$ 
hence we can find $\bar a \in
{}^{\omega >}(M_1),\bar b \in {}^{\omega >}(M_2)$ and $n <
\omega,\varphi(x,\bar b,\bar a)$ such that
\mr
\widestnumber\item{$(iii)$}
\item "{$\circledast(i)$}"  ${\frak C}^{\text{eq}} \models
\varphi[e,\bar b,\bar a]$ 
\sn
\item "{$(ii)$}"   $\models (\exists^{!n} x)\varphi(x,\bar b,\bar a)$ 
\sn
\item "{$(iii)$}"  $n$ minimal under $(i) + (ii)$.
\ermn
We know $\varphi(x,\bar b,\bar a) \vdash \text{ tp}(e,M_1 \cup M_2)$
and let $\{e_0,\dotsc,e_{n-1}\}$ list 
$\varphi({\frak C}^{\text{eq}},\bar b,\bar a)$.

Let $\bar e = \langle e_0,\dotsc,e_{n-1})$.  Possibly increasing $\bar
a,\bar b$ for some formula $\psi = \psi(\bar x,\bar b,\bar a)$ with
$\bar x = \langle x_\ell:\ell <n \rangle$ we have 
$\models \psi(\bar e,\bar b,\bar a)$ and $\psi(\bar x,\bar b,\bar a) 
\vdash \text{ tp}(\bar e,M_1 \cup M_2)$. \nl
So we can find $f$ such that
\mr
\item "{$\circledast$}"  $f$ is an elementary mapping in ${\frak C}$
\nl
Dom$(f) = M_1 \cup M_2 \cup \bar e$ \nl
$f \restriction (M_1 \cup M_2)$ is the identity
\nl
$f(e_0) \ne e_0$ (but of course $f$ permutes $\{e_\ell:\ell < n-1\}$).
\ermn
Let $f(\bar e) = \bar e'$.  Let $\bar e_0 = \bar e,\bar e_1 = f(\bar
e)$. \nl
We can find a sequence of ${\frak C}^{\text{eq}}$-elementary mapping $\langle
g_i:i < |T|^+ \rangle$ such that

$$
\text{Dom}(g_i) = \text{ acl}_{{\frak C}^{\text{eq}}}(M_1 \cup M_2)
$$

$$
g_i \restriction M^{\text{eq}}_2 = \text{ id}
$$

$$
\nonfork{}{}_{M^{\text{eq}}_2} \{\text{Rang}(g_i):i < |T|^+\}. 
$$
\mn
Now
\mr
\item "{$\circledast$}"  if $k < \omega,i_0 < \ldots < i_{k-1} <
\omega$ and $\eta \in {}^n 2$ then the type \nl
$p_\eta = \text{ tp}(g_{i_0}(\bar e_{\eta(0)}) 
\char 94 g_{i_1}(\bar e_{\eta(1)} \char
94 \ldots g_{i_{k-1}}(\bar e_{\eta(k)}),\dbcu_{i < |T|} \text{
Rang}(g_i))$ does not depend on $\eta$.
\ermn
[Why?   By induction on $k$, hence by transitivity of equality it is
enough to prove $p_\eta = p_\nu$ when $1 = |\{\ell:\eta(\ell) \ne
\nu(\ell)\}|$.

By an indiscernible sequence = indiscernible set (= symmetry of nonforking,
etc.) \wilog \, $\eta(0) \ne \nu(0)$.  As Rang$(\bar e_0) = \text{
Rang}(\bar e)$, without loss of generality 
$\dsize \bigwedge_{\ell < k-1} \eta(1 + \ell) = 0 = 
\nu(1 + \ell)$.  Lastly, tp$(\dbcu_{i >0} \text{
Rang}(g_i),\text{Rang}(g_0))$ is finitely satisfiable in $M_2$ so by
the choice of $\psi$ we are done.]

Now for any $\eta \in {}^{(|T|^+)}2$ we define the function $h_\eta$:

$$
\text{Dom}(h_\eta) = M^{\text{eq}}_2 \cup \cup\{g''_i(M^{\text{eq}}_1):
i < |T|^+\} \cup \{g_i(\bar e):i < |T|^+\}
$$

$$
h_\eta \restriction M^{\text{eq}}_2 = \text{ identity}
$$

$$
h_\eta \restriction g''_i(M^{\text{eq}}_1) = \text{ identity}
$$

$$
h_\eta(g_i(\bar e)) = \cases g_i(\bar e) = g_i(\bar e_0) &\text{
\ub{if} } \eta(i) = 0 \\
g_i(\bar e_1) &\text{ \ub{if} } \eta(i)=1 \endcases
$$
\mn
We can find $M_3,M_4,\sigma$ such that

$$
\cup \{g_i(M_1):i < |T|^+\} \subseteq M_3 \prec M_4 \prec {\frak C}
$$

$$
\nonfork{M_2}{M_4}_{M_0}
$$

$$
M_4 \text{ is saturated of cardinality } > \|M_3\|
$$

$$
\sigma \in \text{ Aut}(M_4),\sigma \restriction M_3 = \text{ identity}
$$

$$
(M_4,\sigma) \text{ is a model of } T^{\text{mc}}_{\text{aut}}.
$$
\mn
Now for every $\eta \in {}^{(|T|^+)}2$ we can find $(M^5_\eta,\sigma)
\models T_{\text{aut}}$ such that $(M_4,\sigma) \subseteq (M_5,\sigma)$ and
$\bar b_\eta$ realizing tp$_{{\frak C}^{\text{eq}}}(\bar b,M_0,{\frak C})$
such that

$$
\eta(i) = 0 \Leftrightarrow (\exists \bar x)(\psi(\bar x,\bar
b_\eta,g_i(\bar a)) \wedge \sigma(\bar x) = x).
$$
\mn
So $\{(\exists \bar x)(\psi(\bar x,\bar y,g_i(\bar a)):i < |T|^+\}$ is
an independent set of formulas in $(M_4,\sigma)$ hence
$T^{\text{mc}}_{\text{aut}}$ is unstable.
\bn
\ub{$(a) \Rightarrow (d)$}:

Trivial.
\bn
\ub{$\neg(b) \Rightarrow \neg(c)$}:

Included in the proof of $\neg(b) \Rightarrow \neg(a)$.
\hfill$\square_{\scite{1.1}}$\margincite{1.1}
\enddemo
\bigskip

\demo{\stag{1.2} Observation}  Assume $T^{\text{mc}}_{\text{aut}}$ exists,
$T_*$ any completion of it. \nl
1) If ${\frak C}$ is a saturated model of $T$ of cardinality $\bar
\kappa = \bar \kappa^{< \bar \kappa}$, can be expanded to a model
${\frak C}_*$ of $T_*$. \nl
2) If $M \models T,\sigma \in \text{ Aut}(M)$, let $\sigma^{\text{eq}}$ be the
natural extension of $\sigma$ to an automorphism of $M^{\text{eq}}$, then (it
exists and is unique) $(M^{\text{eq}},\sigma^{\text{eq}}) \models 
(T^{\text{eq}})_{\text{aut}}$. \nl
3) $(T^{\text{eq}})_{\text{aut}}$ 
has a model completion $T$ and there is a natural
one to one correspondence between the completions of the model
completions of $(T^{\text{eq}})_{\text{aut}}$ and $\{T_{**}:T_{**}$ a model
completion of $T^{\text{mc}}_{\text{aut}}\}$ any one of the former is
essentially bi-interpretable with the corresponding one of the latter
(but we have the elements not in any $P_{E(\bar x,\bar y)}$. \nl
4) Let ${\frak C}_* = ({\frak C},\sigma_*)$ be 
a $\bar \kappa$-saturated model of $T_*$
expanding ${\frak C}$.  If $A_\ell \subseteq {\frak C}^{\text{eq}},A_\ell =
\text{ acl}_{{\frak C}^{\text{eq}}}(A_\ell),A_\ell$ 
closed under $\sigma_*,f$ is 
an ${\frak C}^{\text{eq}}$-elementary mapping from $A_1$
onto $A_2$ commuting with $\sigma$ \ub{then} $f$ can be extended to an
automorphism of $({\frak C}^{\text{eq}})_{\text{aut}}$ 
(it is ${\frak C}^{\text{eq}}$ expanded by $\sigma$ 
naturally extended to $\sigma^+$.
\enddemo
\bigskip

\demo{\stag{1.2t} Observation}  1) $M$ is a model of $T,\sigma_* \in
\text{ Aut}(M)$ \ub{iff} $(M,\sigma_*)$ is a model of
$T_{\text{aut}}$. \nl
2) If $M \prec {\frak C}$ and $(M,\sigma_*)$ as a model of
$T_{\text{aut}}$ then for one and only one $\sigma^{\text{eq}}_* \in
\text{ Aut}(M^{\text{eq}})$ extend $\sigma_*$. \nl
3) If $M \prec {\frak C},\sigma^{\text{eq}}_* \in 
\text{ Aut}(M^{\text{eq}})$ then $\sigma^{\text{eq}}_* \restriction M
\in \text{ Aut}(M)$. \nl
4) If $A_\ell \subseteq {\frak C}^{\text{eq}}$ and $A_0 = 
\text{ acl}_{{\frak C}^{\text{eq}}}(A_0)$ and $f_\ell$ is an ${\frak
C}^{\text{eq}}$-elementary mapping from $A_\ell$ onto $A_\ell$ for
$\ell = 0,1,2$  and
$f_0 \subseteq f_1,f_0 \subseteq f_2$ \ub{then} for some automorphism
$F$ of ${\frak C}^{\text{eq}},F \restriction A_0 = \text{ id}_{A_0}$
and $f_2 \cup (F f_1 F^{-1})$ is an elementary mapping in ${\frak
C}^{\text{eq}}$ (hence can be extended to an automorphism of ${\frak
C}^{\text{eq}}$; if $\nonfork{A_1}{A_2}_{A_0}$ then \wilog \, $F
\restriction (A_1 \cup A_2) = \text{ id}_{A_1 \cup A_2}$).
\enddemo
\newpage

\head {\S2} \endhead  \resetall \sectno=2
\bn
\ub{\stag{2.2} Example}:  There is $T$ such that:
\mr
\item "{$(a)$}"  $T$ is as in \scite{1.0}, 
stable $T^{\text{mc}}_{\text{aut}}$ exists.  Moreover $T$ is superstable,
countable \nl
$I(\aleph_\alpha,T) \le 2^{|\alpha|}$ for $\alpha \ge 2^{\aleph_0}$
(hence NDOP, NOTOP, shallow with small depths, with $\le 2^{\aleph_0}$
dimensions)
\sn
\item "{$(b)$}"  $T^{\text{mc}}_{\text{aut}}$ exist
\sn
\item "{$(c)$}"  some completions of $T^{\text{mc}}_{\text{aut}}$ are
stable and some are not.
\endroster
\bigskip

\demo{Proof}  Let us define $M,I$

$|M|$ is $\{(\eta,k,n,\ell):k,n < \omega,\ell < 2$ and $\eta \in
{}^\omega 2\}$ and $k=n \Rightarrow \ell=0$

$E^M_n$, a two-place relation is
$\{(\eta_1,k_1,n_1,\ell_1),(\eta_2,k_2,n_2,\ell_2)) \in |M| \times
|M|:\eta_1 \restriction n = \eta_2 \restriction n\}$

$E^M$, a two-place relation is
$\{(\eta_1,k_1,n_1,\ell_1),(\eta_2,k_2,n_2,\ell_2) \in |M| \times
|M|:\eta_1 = \eta_2\}$

$Q^M$, a one-place relation is
$\{(\eta,k,n,\ell) \in |M|:k = n\}$

$F^M_1$, a one-place relation is:
$F^M_1((\eta,k,n,\ell)) = (\eta,k,k,0)$

$F^M_2$, a one-place relation is:
$F^M_2(\eta,k,n,\ell)) = (\eta,n,n,0)$
\mn
Let $T = \text{ Th}(M)$.  Clearly it satisfies (a):
\mr
\item "{$\circledast_1$}"  $T^{\text{mc}}_{\text{aut}}$ exists. \nl
[Why?  Check that there are no obstructions.]
\sn
\item "{$\circledast_2$}"  $T^{\text{mc}}_{\text{aut}}$ has 
an unstable completion. \nl
[Why?  By \scite{1.1}, or more specifically, see below.]
\ermn
We shall now prove
\mr
\item "{$\circledast^+_2$}"  for $T_*$ a completion of
$T^{\text{mc}}_{\text{aut}},T_*$ is unstable \ub{if}: \nl
for some $M^+ \models T_*$, for some $a \in M^+$ we 
have $\dsize \bigwedge_n a E_n(\sigma^{M^+}(a))$ or just \nl
$(\exists m) \dsize \bigwedge_{n < \omega} a E_n
((\sigma^{M^+})^m(a))$, 
i.e., for some $m^* \in [1,\omega)$ we have $\dsize \bigwedge_n \, a E_n
a_{M^+}$ where $a_0 = a,a_{\ell +1} = \sigma^{M^+}(a_\ell)$ for $\ell <
\omega$.
\ermn
Let $m^*,a,\langle a_\ell:\ell < \omega \rangle$ be as above.
We define $N$ a model of $T$: let $|N|$, the universe of $N$ be

$$
|M^+| \cup \{(m,k,n,\ell):m<m^*,k,n < \omega,\ell <2,k=n \Rightarrow
 \ell=0\}
$$
\mn
we assume no incidental identification.

$$
E^N_n:\cases E^N_n \text{ is an equivalence relation} \\
E^N_n \restriction |M^+| = E^{M^+}_n \\
\text{every } (m,k,n,\ell) \in |N| \backslash |M^+| \text{ is }
E_n\text{-equivalent to } a_m \endcases
$$
\bigskip

$$
E^N:\cases E^N \text{ is an equivalence relation} \\
E^N \restriction |M^+| = E^N \\
\{(m,k,n,\ell) \in |N| \backslash |M^+|:k,n < \omega,\ell < 2,k=n
\Rightarrow \ell=0\} \\
\quad \text{is an } E^N\text{-equivalence class (for each } m < m^*) \endcases
$$

$$
Q^N = Q^N \cup \{(m,k,k,0):k < \omega\}
$$

$$
F^N_1 \text{ extends } F^{M^+}_1,F^N_1((m,k,n,\ell)) = (m,k,k,0)
$$

$$
F^N_2 \text{ extends } F^{M^+},F^N_2((m,k,n,\ell)) = (m,n,n,0).
$$
\mn
Easily
\mr
\item "{$\boxdot_1$}"  $M^+ \restriction \tau_T \prec N$.
\ermn
Now we define an automorphism $\sigma^+$ of $N$:
\mr
\item "{$\boxdot_2$}"  $\sigma^+ \restriction |M^+| =
\sigma^{M^+}$
\sn
\item "{$\boxdot_3$}"  if $m_1,m_2 < m^*,m_2 = m_1 + 1 \text{ mod } m^*$ 
then \nl
$\sigma(m_1,n,k,\ell)$ is: \nl
$(m_2,n,k,1 -\ell)$ if $m_1 = m^* -1 \and n < k$ \nl
$(m_2,n,k,\ell)$ otherwise.
\ermn
Easy to check that $\sigma^+ \in \text{ Aut}(N)$, so $(N,\sigma) \supseteq
M^+$ is a model of $T_{\text{aut}}$.  As
$T^{\text{mc}}_{\text{aut}}$ exists and $M^+ \models
T^{\text{mc}}_{\text{aut}}$ there is a model $N^+ \models
T^{\text{mc}}_{\text{aut}}$ such that $M^+ \prec M^+,(N,\sigma)
\subseteq N^+$.

Let 

$$
\varphi(x,y) = Q(x) \and Q(y) \and x E y \and (\exists z)(F_1(z) 
\and F_2(z) = y \and (\sigma^{m^*}(z) \ne z)]
$$
\mn
This is a first order formula in $\Bbb L(\tau_{\text{Th}(M^+)}) = \Bbb
L(\tau_{T_{\text{aut}}})$ and $N^+ \models \varphi[b_n,b_k]$ iff $n <
\omega$ where $b_n = (0,n,n,0) \in N \subseteq N^+$, so this formula 
has the order property in Th$(N^+) =
\text{ Th}(M^+)$.  So Th$(M^+)$ is unstable as required in
$\circledast^+_2$
\mr
\item "{$\circledast_3$}"  if $T_*$ is a completion of
$T^{\text{mc}}_{\text{aut}}$ not satisfying the demand in $\circledast^+_2$
\ub{then} $T_*$ is stable. \nl
[Why?  As any model $M^+$ of $T_*,\sigma^{M^+}$ acts as a permutation
of $|M^+|/E^{M^+}$ which has no fix point and even no finite cycle.
Now reflect.]
\sn
\item "{$\circledast_4$}"  there is a completion $T_*$ of
$T^{\text{mc}}_{\text{aut}}$ which is stable. \nl
Why?  Let $f$ be a permutation of ${}^\omega 2$ such that
{\roster
\itemitem{ $(\alpha)$ }  $\eta,\nu \in {}^\omega 2 \wedge \eta
\restriction n = \nu \restriction n \Rightarrow f(\eta) \restriction n
= f(\nu) \restriction n$
\sn
\itemitem{ $(\beta)$ }  for every $m<\omega$ ($\ge 2$) for some
$n<\omega$ we have if $\eta \in {}^\omega 2$, then $\eta,f {}^m(\eta)$
are not $E_n$-equivalent.
\endroster}
\ermn
Easy to construct (or use $\dsize \prod_{n<\omega}(n+1)$ instead
${}^\omega 2$) and define $M^+$, a $\tau_{T_{\text{aut}}}$-expansion
of $M$ by defining

$$
\sigma^{M^+}((\eta,k,n,\ell)) = (f(\eta),k,n,\ell)).
$$
\mn
So if $M^+ \subseteq N^+ \models T^{\text{mc}}_{\text{aut}}$ then 
$T_*$ = Th$(N^+)$ fail the demand in $\circledast^+_2$ hence by
$\circledast_3$ it is stable as required 
(and it is uniquely determined by $M^+$,
really just the action on acl$_{{\frak C}^{eq}}(\emptyset)$, suffice.
So $\circledast_4$ holds.  \hfill$\square_{\scite{2.2}}$\margincite{2.2}
\enddemo
\bn
\ub{\stag{2.2t} Discussion}:  It seems reasonable that we can
characterize when this occurs thus answering fully \scite{0.0}; see
below. 
\bn
A closely related example is
\proclaim{\stag{2.3A} Claim}  There is $T$ such that:
\mr
\item "{$(a)$}"  $T$ is stable (complete countable first order theory)
and has elimination of quantifiers for simplicity
\sn
\item "{$(b)$}"  $T$ is superstable and small, i.e., with countable $D(T)$
\sn
\item "{$(c)$}"  $T_{\text{aut}}$ has no model completion
\sn
\item "{$(d)$}"  some $T_{\text{aut}}(M^+)$ has a model completion
where
\endroster
\endproclaim
\bigskip

\definition{\stag{2.3s} Definition}  1) For a model $M^+ = (M,\sigma^{M^+})$ of
$T_{\text{Aut}}$ let $T_{\text{aut}}(M^+) = T_{\text{aut}} 
\cup \text{ Th}(M,c)_{c \in M} \cup \{\sigma(c_1) = 
c_2:\sigma^{M^+}(c_1) = c_1\}$. 
\enddefinition
\bigskip

\remark{\stag{2.3r} Remark}  Actually we 
can use any completion of $T_{\text{aut}} \cup$(the action of
$\sigma$ on acl$_{{\frak C}^{\text{eq}}}(\emptyset,{\frak C}_T)$ (i.e., on the
$E$-equivalence classes for each $n$).
\endremark
\bigskip

\demo{Proof}  Define $M$
\mr
\item "{$(a)$}"  $\tau_M = \{E_n,P_n:n < \omega\} \cup \{E,E_*\}$
\sn
\item "{$(b)$}"  $|M| = \{(\eta,k,n,\ell):\eta \in {}^\omega 2,k <
\omega,n < \omega,\ell < 2\}$
\sn
\item "{$(c)$}"  $E^M_n =
\{(\eta_1,k_1,n_1,\ell_1),(\eta_2,k_2,\eta_2,\ell_2)) \in |M| \times
|M|:\eta_1 \restriction n = \eta_2 \restriction n\}$
\sn
\item "{$(d)$}"  $E^M =
\{(\eta_1,k_1,n_1,\ell_1),(\eta_2,k_2,\eta_2,\ell_2)) \in |M| \times
|M|:\eta_1 = \eta_2$ and $k_1 = k_2\}$
\sn
\item "{$(e)$}"  $E^M_* =
\{((\eta_1,k_1,n_1,\ell_1),(\eta_2,k_2,\eta_2,\ell_2)) \in |M| \times
|M|:\eta_1 = \eta_2,k_1 = k_2,n_1 = n_2\}$
\sn
\item "{$(f)$}"  $P^M_n = \{(\eta,k,n,\ell) \in |M|:n=m\}$.
\ermn
We choose $\sigma^M$ such that $\sigma(\eta,k,n,\ell) =
(\eta',k,n,\ell)$ and $(\eta,\eta')$ are as in the proof of \scite{2.2}.
\enddemo
\bigskip

\remark{Remark}  If we let $(d)'$ be as in \scite{2.5} below we add $\sigma
=$ the identity then $(a) + (c) + (d)'$ is impossibly by
\cite{BlSh:759}.
\endremark
\bn
Actually the case $\sigma$ is the identity on some $M$ is the real one
because
\proclaim{\stag{2.3B} Claim}  For any first order complete $T_1$ (with
$\tau_{T_1}$ a set of predicates for simplicity) there is $T$ such that:
\mr
\item "{$(a)$}"  $T$ is first order complete
\sn
\item "{$(b)$}"  if $a \in M,M \models T$ then we can interpret $T_1$
in $(M,a)$
\sn
\item "{$(c)$}"  $\tau_T \backslash \tau_{T_1}$ countable
\sn
\item "{$(d)$}"  some $T_{\text{aut}}(M^+)$ has a model completion.
\endroster
\endproclaim
\bigskip

\demo{Proof}  As in \scite{2.3A} without $E_*,P_n(n < \omega)$ in any
$E^M$-equivalence class we ``plant" a model of $T_1$.
\enddemo
\bigskip

\proclaim{\stag{2.4} Claim}  Let $T_*$ be a completion of
$T^{\text{mc}}_{\text{aut}}$. \nl
The following are equivalent:
\sn
\ub{Condition $(a)$}:  $T_*$ is stable.
\sn
\ub{Condition $(b)$}:  If $T$ is stable and $(\alpha) + (\beta) +
(\gamma)$ below holds, \ub{then} $(*)$ below holds where
\mr
\item "{$(\alpha)$}"  $M^+_0 \prec M^+_\ell < M^+_3$ for $\ell =
1,2,M_0 \models T_*,M_\ell \models T_{\text{aut}}$ for $\ell = 1,2,3$
and
\sn
\item "{$(\beta)$}"  $M_\ell = M_\ell \restriction \tau_T$ and
$\nonforkin{M_1}{M_2}_{M_0}^{M_3}$ \wilog \, $M_3 \prec {\frak C} =
{\frak C}_T$
\sn
\item "{$(\gamma)$}"  if $f$ is an elementary mapping from 
${\text{\rm acl\/}}_{{\frak C}^{\text{eq}}}(M_1 \cup M_2)$ onto 
itself extending $\sigma^{M^+_1} \cup \sigma^{M^+_2}$ 
\sn
\item "{$(*)$}"  there is an elementary mapping $h$ from 
${\text{\rm acl\/}}_{\frak C},(M_1 \cup M_2)$ onto itself 
such that $h \restriction (M_1 \cup M_2) = { \text{\rm identity\/}}
_{M_1 \cup M_2}$ and $h f h^{-1} =
\sigma^{M^+_3} \restriction { \text{\rm acl\/}}_{{\frak C}^{eq}}
(M_1 \cup M_2)$.
\endroster
\endproclaim
\bigskip

\demo{Proof}  \ub{$(b) \Rightarrow (a)$}:

As in the proof of \scite{1.1}.
\sn
\ub{$\neg(b) \Rightarrow \neg(a)$}:

We can use compactness to replace $\neg(b)$ by a finite failure, and
continue as in the proof of \scite{1.1}.
\enddemo
\bigskip

\remark{\stag{2.5} Remark}  We can make $\neg(b)$ more explicit as in
the proof of \scite{2.4}.
\endremark
\newpage

\head {\S3 NSOP$_3$} \endhead  \resetall \sectno=3
\bigskip

As by \cite{KkSh:748}, if $T^{\text{mc}}_{\text{aut}}$ exists,
\ub{then} $T$ fails the strict order property.  It seems reasonable to
ask if any $T^{\text{mc}}_{\text{aut}}$, which exists, can have the
strict order property.  As we understand the stable case, it seems
reasonable to deal with it.  In fact, more turn out to hold.
\bigskip

\proclaim{\stag{3.1} Claim}  [$T$ as in \scite{1.0}.]  If $T$ is
stable, any completion $T_*$ of $T^{\text{mc}}_{\text{aut}}$ satisfies
NSOP$_3$ (see \cite[\S2]{Sh:500} and \cite{ShUs:789}).
\endproclaim
\bigskip

\demo{Proof}  1) \ub{Clause $(a)$}:

Let $T_*$ be completion of $T^{\text{mc}}_{\text{aut}}$
and $\varphi(\bar x,\bar y)(\ell g(\bar x) = \ell g(\bar y) = n^* <
\omega)$ a first order formula in $\Bbb L(\tau_{T_*})$ and for some $M
\models T_*$ we have $M \models \varphi(\bar a_n,\bar a_m)^{\text{if}(n<m)}$.  Hence we can find an E.M.-template $\Phi$
such that $\tau_\Phi \supseteq \tau_{T_*} = \tau_T \cup \{\sigma\}$
and for linear orders $I \subseteq J$, EM$(I,\Phi) \prec \text{
EM}(J,\Phi) \ne T_*$, with skeleton $\langle \bar a_t:t \in J \rangle$
such that EM$(J,\Phi) \models \varphi[\bar a_s,\bar a_t]^{\text{if}(s <_J
t)}$ for $s,t \in J$ (so $\bar a_t \in \text{ EM}(\{t\},\Phi)$ (see,
e.g., \cite[VII]{Sh:c} or \cite[III]{Sh:e}).  Now (recalling that
EM$_\tau(I,\Phi) = \text{ EM}(I,\Phi) \restriction \tau))$ \wilog \,
\mr
\item "{$\circledast_1$}"  if $I_1,I_2 \subseteq J,I_0 = I_1 \cap I_2$
and if $t \in I_1 \backslash I_0$ then there is $s \in I_0$ such that
$s < t \and ]s,t]_J \cap I_2 \subseteq I_0$ or $t < s \and [t,s]_J
\cap I_2 \subseteq I_0$ then
tp$_{L(\tau_{T^*})}(\text{EM}_{\tau_{T_*}}(I_1,\Phi)$,
EM$_{\tau_{T_*}}(I_2,\Phi))$ is f.s. (finitely satisfiable) in
EM$_{\tau_{T^*}}(I_0,\Phi)$ \nl
[Why?  Let $I \times \Bbb Z$ be ordered lexicographically, choose
$\Phi'$ such that EM$(I,\Phi') = \text{ EM}(I \times \Bbb Z,\Phi)$,
with skeleton $\bar a'_t = \bar a_{(t,0)}$; can look at
\cite{Sh:394}.]
\ermn
For $u \subseteq \{0,1,2\}$ let $M^2_u = \text{ EM}(u,\Phi)$ and if
$|u|=|v|$ both subsets of $\{0,1,2\}$ let $f_{v,u}$ be the canonical
isomorphism from $M_u$ onto $M_v$.  Let $M^1_u = M^2_u \restriction
\tau_{T_*},M^0_u = M^2_u \restriction \tau_T$.  Let $N$ be such that
$M^0_{\{0,1,2,\}} \prec N,N$ is $\|M^0_{\{0,1,2,\}}\|^+$-saturated
\mr
\item "{$\circledast_2$}"  in
$N,\nonfork{}{}_{M^0_\emptyset}\{M^0_{\{0\}},M^0_{\{1\}},M^0_{\{2\}}\}$
\nl
[Why?  By $\circledast_1$ and nonforking calculus.]
\ermn
Let $g_0 =: f_{\{0\},\{2\}} \cup f_{\{2\},\{0\}}$
\mr
\item "{$\circledast_3$}"  $g_0$ is an elementary mapping (inside $N$)
\nl
[Why?  Nonforking calculus.]
\ermn
Let $g_1$ be an elementary mapping inside $N$ extending $g_0$ with
domain $M^0_{\{0,2\}}$.

Let $M^{0,*}_{\{0,2\}} = g(M^0_{\{0,1\}})$.

Let $M^{1,*}_{\{0,2\}}$ be an expansion of $M^{0,*}_{\{0,2\}}$ by an
automorphism $\sigma^{M^{1,*}_{\{0,2\}}}$ such that $g_1$ is an
isomorphism from $M^1_{\{0,2\}}$ onto $M^{1,*}_{\{0,2\}}$, clearly
exists.

As $N$ is a model of the stable theory $T$ without loss of generality \nl
tp$_{L^*(\tau_T)}(|M^{1,*}_{\{0,2\}}|,|M^0_{\{0,1,2\}}|)$ does not
fork over $|M^0_{\{0\}}| \cup |M^0_{\{2\}}|$.

Now the point is that
\mr
\item "{$\bigodot$}"  $h = \sigma^{M^1_{\{0,1\}}} \cup
\sigma^{M^{1,*}_{\{0,2\}}} \cup \sigma^{M^1_{\{1,2\}}}$ is a
permutation of $|M^{1,*}_{\{0,1\}}| \cup |M^1_{\{0,1\}}| \cup
|M^1_{\{1,2\}}|$ and is an elementary mapping. \nl
[Why?  Let $B_0 = |M^0_{\{0\}}| \cup |M^0_{\{2\}}|,B_1 =
|M^0_{\{0,1\}}| \cup |M^0_{\{2,2\}}|$.
\nl
By \cite[XII]{Sh:c}, the pair $(B_0,B_1)$ satisfies the T.V. condition
inside $N$ (i.e., if $\varphi(\bar x,\bar y) \in \Bbb L(\tau_T),N
\models \varphi[\bar a,\bar b],\bar a \subseteq B_1,\bar b \subseteq
B_0$ then for some $\bar a' \subseteq B_0,N \models \varphi[\bar
a',\bar b]$.  Moreover, we can allow $\bar b \subseteq
|M^{0,*}_{\{0,2\}}|$ then this follows.] 
\ermn
So for some $N',N \prec N' \models T$ and there is an automorphism $h'$ of
$N'$ extending $h$ and we can extend $(N',h')$ to a model $(N'',h'')$
of $T_*$.  By this model clearly

$$
(N'',h'') \models \varphi[\bar a_0,a_1] \text{ using } M^1_{\{0,1\}}
$$

$$
(N'',h'') \models \varphi[\bar a_1,\bar a_2] \text{ using } M^1_{\{1,2\}}
$$

$$
\align
(N'',h'') \models &\varphi[\bar a_2,\bar a_0] \text{ using } 
M^{1,*}_{\{0,2\}} \text{ and} \\
  &g_1 \text{ being an isomorphism from } M^1_{\{0,2\}} \text{ onto }
M^{1,*}_{\{0,2\}}.
\endalign
$$
\mn
This is enough to show $T_* \models \text{ NDOP}_3$.
\enddemo
\bigskip

\proclaim{\stag{3.2} Claim}  $T$ is stable or just simple \ub{then} 
any $T_*$ (assuming it exists, $K_*$ in general) is simple.
\endproclaim
\bigskip

\demo{Proof}  We write it for $K_*$.
Choose $\kappa = \text{ cf}(\kappa) > |T|$ and $\mu$ a
strong limit singular cardinal of cofinality $\kappa$.  Let $\langle
\lambda_i:i < \kappa \rangle$ be increasing with limit $\mu,\lambda_0
> \kappa,\lambda_\kappa = \mu,\langle {}_*M^+_i:i < \kappa \rangle$ be
an increasing sequence of elementary submodels of ${\frak C}_{K_*}$
(check notation), $\|{}_*M^+_i\| = 2^{\lambda_i},{}_*M^+_i$ is
$\lambda^+_i$-homo universal (in $K^{\text{ec}}_{\text{aut}}(T)$),
$M^+ = \cup\{{}_*M^+_i:i < \kappa\}$.  Let $\langle p^+_i:i < \mu^+
\rangle$ be a sequence of existential types in $\Bbb L(\tau
\cup\{\sigma\})$ each of cardinality $\le \kappa$ with domain
$\subseteq M$, and we shall prove that for some $\alpha < \beta <
\mu^+,p^+_\alpha \cup p^+_\beta$ is realized in ${\frak C}_{K_*}$,
this suffices.

For each $\alpha < \mu^+$, we can find $a_\alpha \in {\frak C}_{K_*}$
realizing $p_i$ and $N^*_{3,\alpha} \prec {\frak C}_{K_*}$ of
cardinality $\kappa$ to which $a_i$ belongs and $N^+_{2,\alpha} =
N^+_{3,\alpha} \cap M^+ \prec M^+$ and tp$_{\frak
C}(|N^+_{3,\alpha}|,|M^+|)$ does not fork over $|N^+_{2,\alpha}|$.
Let $N^+_{1,\alpha} \prec N^+_{3,\alpha}$ be of cardinality $|T|$ such
that $a_i \in N^+_{1,\alpha}$, tp$_{\frak C}(|N^+_{1,\alpha}|,|M^+|)$
does not fork over $|N^+_{0,i}$ where $N^+_{0,i} = N^+_{1,i}
\restriction M^+ \prec M^+$.  Without loss of generality $\alpha <
\mu^+ \Rightarrow N^+_{0,\alpha} = N^+_0$ and for every $\alpha,\beta
< \mu^+$ there is an isomorphism $h_{\beta,\alpha}$ from
$N^+_{3,\alpha}$ onto $N^+_{3,\beta}$ mapping
$a_\alpha,N^+_{1,\alpha},N^+_{2,\alpha}$ to
$a_\beta,N^+_{1,\beta},N^+_{2,\beta}$ respectively and
$h_{\beta,\alpha} \restriction N^+_0 = \text{ id}_{N^+_0}$.  Moreover,
\wilog \, for some well ordering $<^*$ all $h_{\beta,\alpha}$ are
order preserving. \nl
Let $\kappa > \bar \kappa,{\frak B}$ 
be an elementary submodel of $({\Cal H}(\chi),\in)$
of cardinality $2^\kappa$ such that $T,\kappa,\mu,{\frak C},{\frak
C}_{K_*},M^+,\langle N^+_i:i < \mu^+ \rangle$ belongs and such that
$[{\frak B}] \le \kappa \subseteq {\frak B}$.  Now choose $\alpha(2)
\in \mu^+ \backslash {\frak B}$, and let $M^+_0 = N^+_{1,\alpha}
\restriction {\frak B}$.  Clearly $M^+_0 \prec M^+$ and there is
$\alpha(1) \in \mu^+ \cap {\frak B}$ such that
$h_{\alpha(1),\alpha(2)}$ is the identity on $M^+_0$. \nl
[FILL?]
\enddemo
\newpage

     \shlhetal 

\newpage
    
REFERENCES.  
\bibliographystyle{lit-plain}
\bibliography{lista,listb,listx,listf,liste}

\enddocument